\documentclass[12pt]{article}

\usepackage{xcolor}
\usepackage{amssymb}
\usepackage{amsfonts}
\usepackage{amsmath}
\usepackage{enumerate}
\usepackage{tabu}
\usepackage{yfonts}
\usepackage{mathrsfs}
\usepackage{bbm}
\usepackage[margin=1in]{geometry}

\newcommand\F{\mathbb{F}}
\newcommand\Q{\mathbb{Q}}

\newcommand\R{\mathbb{R}}

\newcommand\cB{\mathcal{B}}
\newcommand\cC{\mathcal{C}}
\newcommand\cF{\mathcal{F}}
\newcommand\cG{\mathcal{G}}
\newcommand\cH{\mathcal{H}}
\newcommand\cJ{\mathcal{J}}
\newcommand\cM{\mathcal{M}}

\newcommand\al{\alpha}
\newcommand\bt{\beta}

\newcommand\Gm{\Gamma}

\newcommand\eps{\epsilon}
\newcommand\lm{\lambda}
\newcommand\sg{\sigma}

\newcommand\Dl{\Delta}
\newcommand\Lm{\Lambda}

\newcommand\la{\langle}
\newcommand\ra{\rangle}
\newcommand\lla{\la\!\la}
\newcommand\rra{\ra\!\ra}
\newcommand\ad{\mathrm{ad}}

\newcommand\Aut{\mathrm{Aut}}
\newcommand\Miy{\mathrm{Miy}}
\newcommand\End{\mathrm{End}}

\newtheorem{thm}{Theorem}[section]

\newtheorem{defn}[thm]{Definition}
\newtheorem{exa}[thm]{Example}

\newtheorem{prob}[thm]{Problem}

\title{Axial Algebras: Questions and Conjectures}
\author{I.~Gorshkov, S.~Shpectorov}

\date{\vspace{-30px}}

\begin{document}
\maketitle

\begin{abstract}
Axial algebras are non-associative algebras generated by idempotents, called axes, whose adjoint action 
satisfies a fusion law. When this fusion law is graded, axes naturally lead to automorphisms of the algebra, 
and so such axial algebras are inextricably linked with groups. 

This article is meant to complement the recent survey \cite{ms} by significantly expanding the list of 
interesting open problems suggested by the specialists in the field, and providing a further discussion 
of the related concepts and available results.
\end{abstract}
\newcommand{\Addresses}{{
		\bigskip\noindent
		\footnotesize
		Ilya~Gorshkov, \textsc{Sobolev Institute of Mathematics, Novosibirsk, Russia;}\\\nopagebreak
		\textit{E-mail address: } \texttt{ilygor8@gmail.com}
		
		\medskip\noindent
		Sergey~Shpectorov, \textsc{School of Mathematics, University of Birmingham,
		Edgbaston, Birmingham B15 2TT, United Kingdom;}
		\\\nopagebreak
		\textit{E-mail address: } \texttt{s.shpectorov@bham.ac.uk}
    	\medskip
}}

\section{Introduction}

Algebras considered here are non-associative, which means that associativity of the product is not required. 
Hence an algebra over a ring of coefficients $k$ is a $k$-module endowed with a $k$-bilinear product 
operation. In most cases, $k$ is a field; however, the more general 
rings of coefficients also arise in specific situations.

The axial paradigm is a novel way to define interesting classes of algebras. Normally, these are defined in 
terms of identities that must be satisfied in the target class of algebras. Under the axial approach, we 
consider algebras generated by non-zero idempotents, called axes, whose adjoint action is governed a fusion 
law. Hence the fusion laws can be used to delineate interesting classes of algebras. Formally, the axial 
approach was introduced in \cite{hrs1,hrs2}, however, as we will soon see, axial behaviour has been observed 
in algebras since a long time ago.

Development of the theory of algebras has always been aligned with the study of groups. Hence it is no 
surprise that axial algebras are also closely related to groups. The motivating example was the Griess algebra 
for the Monster finite sporadic simple group $M$. Recall that $M$ was constructed by Griess \cite{g} as the 
automorphism group of a commutative non-associative algebra $V$ of dimension $196884$ over $\R$, which is now known as the Griess algebra. This algebra $V$
appears as weight-2 component of the Moonshine vertex operator algebra (VOA) $V^\natural$ which is one of the 
central objects in the theory of VOAs. Miyamoto \cite{m} discovered that the $2A$ involution from $M$ arises in 
a natural way from the conformal Ising vectors contained in $V$. These Ising vectors are the same, up to a scalar, as the $2A$-axes, found by Conway and Norton within the context of the Monster Moonshine 
conjecture. 

Miyamoto \cite{m,m03} started the programme of classifying VOA generated by two Ising vectors. This was completed 
by Sakuma \cite{s} who managed to obtain a classification of the weight-2 component algebras in such VOA. 
As it turned out all of them arise as subalgebras of the Griess algebra generated by two $2A$-axes. Following this, 
Ivanov \cite{i} realised that the result of Sakuma can be reformulated as a theorem on the class of algebras whose axioms are modelled on the properties of the Griess algebra $V$. He formulated the suitable axioms and hence introduced the class of Majorana algebras, which are commutative non-associative algebras 
generated by special idempotents called Majorana axes. Within these axioms, the remarkable Table \ref{Monster432} appeared.
\begin{figure}[ht]
\begin{center}
\renewcommand{\arraystretch}{1.6}
\begin{tabular}{|c||c|c|c|c|}
  \hline
   ${*}$ & $1$ & $0$ & $\frac{1}{4}$ & $\frac{1}{32}$ \\
   \hline \hline
  $1$ & $1$ &  & $\frac{1}{4}$ & $\frac{1}{32}$ \\
  \hline
  $0$ &  & $0$ & $\frac{1}{4}$& $\frac{1}{32}$ \\
  \hline
  $\frac{1}{4}$ & $\frac{1}{4}$ &$\frac{1}{4}$ & $1,0$& $\frac{1}{32}$  \\
   \hline
   $\frac{1}{32}$ & $\frac{1}{32}$ & $\frac{1}{32}$ & $\frac{1}{32}$ & $1,0,\frac{1}{4}$\\
   \hline
\end{tabular}
\caption{Monster fusion law}\label{Monster432}
\end{center}
\end{figure}

This table represents what is now called a fusion law, specifically, the Monster fusion law, and it describes the properties of the adjoint action of Majorana axes. Up to our knowledge, this was the first time a fusion law appeared as part of the axiomatics of a class of algebras. However, similar properties of idempotents in algebras have been observed earlier. 

In 1881, Peirce \cite{p} described how an associative algebra decomposes with respect to an idempotent. He considered the general non-commutative case, where the left and right adjoint actions of the idempotent are different, and his decomposition consisted of four pieces. In the commutative case, this simplifies to a two-piece decomposition described by Table \ref{Associative law}.
\begin{figure}[ht] 
\begin{center}
\renewcommand{\arraystretch}{1.3}
\begin{tabular}{|c||c|c|}
  \hline
   ${*}$ & $1$ & $0$ \\
   \hline \hline
  $1$ & $1$ & \\
  \hline
  $0$ &  & $0$ \\
  \hline
\end{tabular}
\caption{Associative fusion law}\label{Associative law}
\end{center}
\end{figure}

In 1947, Albert introduced three-term Peirce decompositions in Jordan algebras. Again, these are with respect to an arbitrary idempotent and they can be described by Table \ref{Jordan fusion law}.
\begin{figure}[ht]
\begin{center}
\renewcommand{\arraystretch}{1.8}
\begin{tabular}{|c||c|c|c|c|}
  \hline
   ${*}$ & $1$ & $0$ & $\frac{1}{2}$ \\
   \hline \hline
  $1$ & $1$ &  & $\frac{1}{2}$ \\
  \hline
  $0$ &  & $0$ & $\frac{1}{2}$ \\
  \hline
  $\frac{1}{2}$ & $\frac{1}{2}$ & $\frac{1}{2}$ & $1,0$\\
  \hline
\end{tabular}
\caption{Jordan fusion law}\label{Jordan fusion law}
\end{center}
\end{figure}
Note that the actual tables, as above, do not appear in the works of Peirce and Albert. Instead, the properties represented by the table are given by formulas in the text. The table simply expresses the same properties in a more compact form.

In the background chapter we will introduce the basic concepts and definitions. In cases where we need more specific definitions to formulate a question, we will provide them in the text.

All questions are divided into chapters according to the specifics of the question.

\section{Background}
This section introduces the foundational definitions of axial algebras and their subclasses, focusing on algebras of Jordan and Monster type.

\begin{defn}
A fusion law is a pair $(\cF,\ast)$, where $\cF$ is a set and $\ast:\cF\times\cF\to 
2^\cF$ is a binary operation on $\cF$ with values in the set $2^\cF$ of all subsets of 
$\cF$. 
\end{defn} 

The fusion laws we focus on are defined on a small set $\mathcal{F}$, and we find it useful to organise them in a table similar to a group multiplication table. To simplify the notation within the table, we will omit the set brackets and simply list the members of the resulting set $\lambda \star \mu$ directly within the $(\lambda, \mu)$ entry of the table. Specifically, if the entry is left blank, it indicates that the set $\lm\ast\mu$ is empty.

Figure \ref{Jordan type law} shows the fusion law of Jordan type $\eta$, denoted by $\cJ(\eta)$. Here $\cF=\cJ(\eta)=\{1,0,\eta\}$, $\eta\neq 1,0$, is a subset of a field $\F$. The fusion law $\cJ(\eta)$ generalises the Jordan fusion law from Figure \ref{Jordan fusion law}, which arises when $\eta=\frac{1}{2}$.

\begin{figure}[ht] 
\begin{center}
    \renewcommand{\arraystretch}{1.5}
	\begin{tabular}{|c||c|c|c|}
		\hline
		$\ast$ & $1$ & $0$ & $\eta$\\
		\hline\hline
		$1$ & $1$ & & $\eta$\\
		\hline
		$0$ & & $0$ & $\eta$\\
		\hline
		$\eta$ & $\eta$ & $\eta$ & $1,0$\\
		\hline
	\end{tabular}
	\caption{Jordan type fusion law $\cJ(\eta)$}\label{Jordan type law}
\end{center}
\end{figure}

Another example of a fusion law is the Monster type fusion law $M(\alpha,\beta)$, see Figure \ref{Monster law}. Here $\al$ and $\bt$ are distinct elements of $\F$ different from $0$ and $1$. When $\alpha=\frac{1}{4}$ and $\beta=\frac{1}{32}$, the Monster type fusion law $M(\alpha,\beta)$ specialises to the Monster fusion law (see Figure \ref{Monster432}).
\begin{figure}[ht] 
\begin{center}
    \renewcommand{\arraystretch}{1.4}
\begin{tabular}{|c||c|c|c|c|}
  \hline
   $\ast$ & $0$ & $1$ & $\al$ & $\bt$ \\
   \hline \hline
  $0$ & $0$ & & $\al$ & $\bt$ \\
  \hline
  $1$ & & $1$ & $\al$& $\bt$ \\
  \hline
  $\al$ & $\al$ &$\al$ & $1,0$& $\bt$  \\
   \hline
   $\bt$ & $\bt$ & $\bt$ & $\bt$ & $1,0,\al$\\
   \hline
\end{tabular}
\caption{Monster type fusion law $\cM(\al,\bt)$}\label{Monster law}
\end{center}
\end{figure}

Let $A$ be a commutative (non-associative) algebra over a ﬁeld $\F$.

\begin{defn}
\begin{enumerate}
\item[(a)] For $a\in A$, the (left) adjoint action $ad_a : A \rightarrow A$ is given by $ad_a(b)=ab$ for all $b\in A$.
\item[(b)] For $\lm\in\mathbb{F}$, set $A_{\lambda}(a)=\{b \in A | ad_a(b) =\lambda b\} = \{b \in A | ab = \lambda b\}$.
\item[(c)] For $\Lambda\subseteq\mathbb{F}$, set $A_{\Lambda}(a) = \bigoplus_{\lambda\in\Lambda}A_{\lambda}(a)$.
\end{enumerate}
\end{defn}

Let $\mathcal{F}\subseteq\mathbb{F}$ be a fusion law.
 
\begin{defn}
An idempotent $0\neq a\in A$ is an \emph{$\mathcal{F}$-axis} if 
\begin{enumerate}
\item[(a)] $A = A_{\mathcal{F}}(a)$; and 
\item[(b)] $A_{\lm}(a)A_{\mu}(a)\subseteq A_{\lm\ast\mu}(a)$ for all $\lm, \mu \in \cF$. 
\end{enumerate}
\end{defn}

If it is clear which fusion law we are considering, we will simply call $a$ an axis.

This definition assumes that $1\in\cF$, as $ad_a(a) = aa = a = 1a$. 

\begin{defn}
An axis a is \emph{primitive} if $A_1(a) = \la a\ra.$
\end{defn}

\begin{defn}
An algebra $A$ over $\F$ is called (primitive) \emph{$\cF$-axial} if it is commutative and generated by (primitive) $\cF$-axes.
\end{defn}

Often we view an axial algebra as a pair $(A,X)$, where $X$ is a specific set of (primitive) axes generating $A$.

We will mainly consider primitive axes and axial algebras, so for simplicity we will usually talk about axes and axial algebras, assuming their primitivity. Where we encounter non-primitive axes, this will be explicitly stated.

For a subset $X\subseteq A$, we denote the subalgebra generated by $X$ by $\lla X\rra$. We will say that an axial algebra $A$ is $n$-generated if $A$ is generated by $n$ axes.

Often algebras admit bilinear forms that associate with the algebra product.

\begin{defn} \label{Frobenius}
Suppose $A$ is an axial algebra. A non-zero bilinear form $(\cdot, \cdot)$ on $A$ is called a \emph{Frobenius form} if
$$(uv, w) = (u, vw)$$
for all $u, v, w \in A$.
\end{defn}

The term ``Frobenius form'' was borrowed from the theory of Frobenius algebras. These are associative algebras equipped with a non-degenerate bilinear form associating with the algebra product, as in the above definition. However, it is important to note that in Definition \ref{Frobenius} we do not require the Frobenius forms to be non-degenerate. Indeed, the radical of the Frobenius form is important for the structural theory of axial algebras.

We now turn to graded fusion law. In the known examples, we mostly have grading by the group $C_2$ of order $2$, but a general abelian grading group $T$ can be introduced in much the same way.

\begin{defn}
Suppose that $\cF$ is a fusion law and $T$ is an abelian group. A \emph{$T$-grading} of $\cF$ is a partition $\{ \cF_t : t \in T \}$ of $\cF$, with some parts $\cF_t$ possibly empty, such that $\cF_s \ast \cF_t \subseteq \cF_{st}$ for all $s, t \in T$. 
\end{defn}

When $T=\la t\ra\cong C_2$, instead of $\cF_t$ and $\cF_1$ we write $\cF^-$ and $\cF^+$, respectively. Note that the law $\cJ(\eta)$ is $C_2$-graded with $\cJ(\eta)^-=\{\eta\}$ and $\cJ(\eta)^+=\{0,1\}$. Also, the Monster type fusion law $\cM(\alpha,\beta)$ is $C_2$-graded with $\cM(\alpha,\beta)^-=\{\beta\}$ and $\cM(\alpha,\beta)^+=\{0,1,\alpha\}$.

Given a $T$-grading of $\cF$ and an axis $a$ in an $\cF$-axial algebra $A$, it is immediate that
$$A = \bigoplus_{t \in T} A_{\cF_t}(a)$$
is a $T$-grading of the algebra $A$. It now follows that, for each linear character $\chi$ of $T$, we can define an automorphism of $A$. Namely, let $\tau_a(\chi)$ be the linear map $A \to A$, which acts on $A_t(a) := A_{\cF_t}(a)$ by multiplying with the scalar $\chi(t)$. The map $\tau_a(\chi)$ is an automorphism of $A$ called a \emph{Miyamoto automorphism}. E.g., when $T=\la t\ra\cong C_2$, we have a unique non-trivial linear character $\chi$ and $\chi(t)=-1$. Then $\tau_a(\chi)$ negates $A_-=A_{\cF_-}(a)$ and acts as identity on $A_+=A_{\cF_+}(a)$. For brevity, instead of $\tau_a(\chi)$ we will then write simply $\tau_a$.

It is immediate that the map $\chi \mapsto \tau_a(\chi)$ is a homomorphism from the group $T^*$ of linear characters of $T$ to $\operatorname{Aut}(A)$. The image of this map, $T_a \leq \operatorname{Aut}(A)$, is called the axial subgroup corresponding to the axis $a$. 

\begin{defn}
For a $T$-graded fusion law $\cF$ and an $\cF$-axial algebra $A$ with the set of generating axes $X$, the group
$$\Miy(A, X) = \la T_a : a \in X \ra = \la \tau_a(\chi) : a \in X, \chi \in T^* \ra \leq \Aut(A)$$
is called the \emph{Miyamoto group} of $A$. We will also use the notation $\operatorname{Miy}(A)$ and $\operatorname{Miy}(X)$, depending on what we focus on.
\end{defn}

\subsection{3-Transposition groups and Matsuo algebras}

In this subsection we introduce the class of Matsuo algebras which form the bulk of algebras of Jordan type. However, we need to start with the discussion of $3$-transposition groups 

\begin{defn}
An $n$-transposition group $(G,D)$ is a group $G$ generated by a normal set 
$D$ of involutions such that: $|de|\leq n$ for all $d,e \in D$. 
\end{defn}

The most well-studied case is the case of $n=3$, that is, the case of $3$-transposition groups. The motivating example of such groups is given by the symmetric group $G=S_n$ generated by its class of \emph{transpositions} 
($2$-cycles) $D=(1,2)^G$. Two transpositions $\sg$ and $\tau$ are either equal, in which case $|\sg\tau|=1$, or their supports are disjoint, in which case they commute, that is $|\sg\tau|=2$, or their support sets meet in one element, in which case they do not commute and $|\sg\tau|=3$.

The class of $3$-transposition groups was first studied by Fischer, who classified them in an important special case \cite{f}. To formulate his result we will need an additional definition.

\begin{defn}
Given a $3$-transposition group $(G,D)$ and a set $C\subset D$ of $3$-transpositions, the \emph{diagram} of $C$ is the graph $(C)$ on the vertex set $C$, where two vertices $c,d\in C$ are connected by an edge if and only if $|cd|=3$.
\end{defn}

When the diagram $(D)$ of the full set $D$ is disconnected, say with connected components $D_i$, $i\in I$, we have that $G$ is the central product of the smaller $3$-transposition groups $(G_i,D_i)$, where $D_i:=\la D_i\ra$. In view of this, we are mostly interested in the \emph{connected} $3$-transposition groups, which are the groups $(G,D)$ with the property that the diagram $(D)$ is a connected graph.

Let us also mention that it is easy to see that the connected components $D_i$ of $D$ are exactly the conjugacy classes of $G$ contained in the normal set $D$. Hence the connected $3$-transposition groups are exactly those where $D$ is a single conjugacy class of $G$.

We will now state the classification of $3$-transposition groups in the important special case, which is generally attributed to Fischer, even though his published result in \cite{f2} is slightly weaker.

\begin{thm}
Suppose that $(G,D)$ is a finite connected $3$-transposition group such that every soluble normal subgroup of $G$ is contained in the centre $Z(G)$ of $G$. 
Then $G/Z(G)=1$ (and so $G\cong C_2$) or it is one of the following groups:
\begin{enumerate}
\item[(a)] a symmetric group $S_n$ of degree $n$ for $n \geq 5$;
\item[(b)] a symplectic group $\operatorname{Sp}(2n, 2)$ over a field with $2$ elements for $n \geq 2$;
\item[(c)] an orthogonal group $\operatorname{O}^{\mu}(2n, 2)$ for $\mu \in \{-1, 1\}$  and $n \geq 4$;
\item[(d)] a unitary group $\operatorname{PSU}(n, 2)$ over a field with $4$ elements for $n \geq 4$;
\item[(e)] an orthogonal group $\operatorname{O}^\mu(n, 3)$ over a field with $3$ elements for $\mu\in\{-1, 1\}$ and $n \geq 4$; or
\item[(f)] one of the five groups $\operatorname{O}^+(8,2):S_3$, $\operatorname{O}^+(8,3):S_3$, $Fi_{22}$, $Fi_{23}$, $Fi_{24}$.
\end{enumerate}
\end{thm}

The general classification of $3$-transposition groups was achieved by Cuypers and Hall \cite{hc}. The general statement is too complicated to include here. However, the important point is that it includes the reducible case. We will now explain what this means.

The general $3$-transposition group $G$ leaves invariant two equivalence relations on the set $D$. For $c\in D$, we write $A_c$ for the set of all $d\in D$ such that $|cd|=3$. Similarly, we write $D_c$ for the set of all $d\in D$ such that $|cd|=2$. Then we have the equivalence relation $\tau$ and $\theta$ defined by: $c\tau d$ if and only if $A_c=A_d$ and $c\theta d$ if and only if $D_c=D_d$. When both these equivalences are the trivial ones, we say that the $3$-transposition group $(G,D)$ is irreducible, and in fact, the condition in the theorem above means exactly that $(G,D)$ is not only connected but also irreducible. Reducible $3$-transposition groups have normal soluble subgroups non contained in the centre. Namely, if $\tau$ is nontrivial then 
$$\tau(G)=\la cd\mid c,d\in D,c\tau d\ra$$
is a nontrivial and non-central normal $2$-subgroup. Symmetrically, if $\theta$ is nontrivial then
$$\theta(G)=\la cd\mid c,d\in D,c\theta d\ra$$
is a nontrivial and non-central normal $3$-subgroup of $G$. Note that only one of these two normal subgroups can be non-trivial, that is, only one of the two equivalences $\tau$ and $\theta$ can be nontrivial. 

Factoring out $\tau(G)$ or $\theta(G)$ (whichever one of them is nontrivial), we obtain a smaller $3$-transposition group, and this is why we call such a $3$-transposition group reducible. Note that the group obtained by such a reduction can still be reducible, and so additional reductions can be required until we reach an irreducible $3$-transposition group from Fischer's list.

\medskip
Let us now describe the important class of axial algebras related to $3$-transposition groups.

\begin{defn}\label{Matsuo}
Let $(G, D)$ be a group of $3$-transpositions. 
Suppose that $\F$ is a field of characteristic not $2$ and $\eta \in \F \setminus \{0, 1\}$. The \emph{Matsuo algebra}  
$A = M_{\eta}(G, D)$ has basis $D$ and multiplication given by  
$$
a \circ b = 
\begin{cases} 
a, & \text{if } b = a, \\
0, & \text{if } |ab| = 2, \\
\frac{\eta}{2}(a + b - c), & \text{if } |ab| = 3.
\end{cases}
$$
Here $a,b\in D$ and, when $|ab|=3$, $c=a^b=b^a$.
\end{defn}

Clearly, the elements of $D$, forming the basis of $M=M_\eta(G,D)$ 
are idempotents. It is easy to see that they are primitive and they 
satisfy the fusion law $\cJ(\eta)$ (see Figure \ref{Jordan type law}). 
Hence, Matsuo algebras belong to the class of algebras of Jordan type 
$\eta$. Two specific examples of Matsuo algebras that will appear below are 
the algebra $2B=\F\oplus\F$ constructed from a pair of computing 
involutions in a group $G$ of order $4$, and the algebra $3C(\eta)$ 
constructed from the class of three involutions in $G=S_3$.

Every Matsuo algebra admits a Frobenius form, which is given by
$$
(a,b)= 
\begin{cases} 
1, & \text{if } b = a, \\
0, & \text{if } |ab| = 2, \\
\frac{\eta}{2}, & \text{if } |ab| = 3.
\end{cases}
$$

Before we leave the topic of $3$-transposition groups and Matsuo algebras, let us mention the geometric concepts related to them.

\begin{defn}
A Fischer space is a point-line geometry with all lines of size $3$, which is a partial linear space (i.e., every pair of points lines in at most one line), and finally, such that every pair of intersecting lines generates either the dual affine plane of order $2$ (the dual of the Fano plane) or the affine plane of order $3$.
\end{defn}

This concept was introduced by Buekenhout \cite{bu}, who observed that such a geometry can be constructed from every $3$-transposition group $(G,D)$. Namely, we take $D$ as the set of points and all triples $\{a,b,c\}$ with $|ab|=3$ and $c=a^b=b^a$ as lines. In fact, every Fischer space arises in this way from a $3$-transposition group. Furthermore, the group $(G,D)$ can be recovered from its Fischer space up to centre, that is, Fischer spaces bijectively correspond to $3$-transposition groups taken up to central equivalence.

Fischer spaces is a very interesting class of geometries and in fact this geometric language was instrumental in the classification of $3$-transposition groups by Cuypers and Hall \cite{hc}. Finally, the Matsuo algebra can be defined in terms of the Fischer space, which means that centrally equivalent $3$-transposition groups produce isomorphic Matsuo algebras.

\subsection{Baric algebras}

In almost every class of axial algebras we have special examples admitting an additional structure. An algebra $A$ over $\F$ is said to be \emph{baric} if it admits a surjective algebra homomorphism $w:A\to\F$. This concept was introduced by Etherington \cite{e} within the context of genetic algebras. The homomorphism $w$ is usually called the \emph{baric weight} (or sometimes, the augmentation function) and its kernel is called the baric ideal or \emph{baric radical} (also, augmentation ideal).
Clearly, the baric radical has codimension $1$ and it is indeed an ideal of $A$.

Note that $A$ may admit many different baric weights. Things change when we step into the axial world.

\begin{defn}
An axial algebra $(A,X)$ is said to be \emph{baric} if it admits a baric weight $w:A\to\F$ such that $w(a)=1$ for all generating axes $a\in X$.
\end{defn}

Since $a\in X$ is an idempotent, $w(a)$ is an idempotent in $\F$, which means that $w(a)=0$ or $1$. So the extra condition above simply means that $w(a)\neq 0$ for all $a\in X$, that is, that $X$ meets the baric ideal trivially.

It is easy to see that such a baric structure on an axial algebra is unique when it exists. 

Every baric axial algebra admits a Frobenius form. Indeed, we define it as
$(u,v):=w(uv)$
for $u,v\in A$. It is easy to see that this is a symmetric bilinear form. Furthermore, $(a,a)=1$ for each $a\in X$. Also, the radical of this Frobenius form coincides with the baric radical.

\section{General axial algebras}

\subsection{Frobenius form}

As we have already mentioned, a great majority of known axial algebras admit a Frobenius 
form.

\begin{prob}   
Does a Frobenius form exist on a general axial algebra? When is it non-singular? positive 
(non-negative) definite?
\begin{flushright} V.G. Tkachev\end{flushright}
\end{prob}

Of course, the final part of this question makes senseonly  when we work over the field of real numbers or one of its subfields.

While it is difficult to expect that every axial algebra admits a (non-zero) Frobenius form, 
it may be the case for specific classes of axial algebras. For example, in \cite{hss1} it 
was shown that every algebra of Jordan type admits a unique Frobenius form such that every 
primitive axis has length (Frobenius norm) $1$.

Sometimes existence of a Frobenius form is included in the axiomatics of the specific class 
of axial algebras. Also, axial algebras arising in applications often come with a Frobenius 
form, which carries important information about the algebra.

\begin{prob}
Normally, the (primitive) axes are a small part of the complete set of idempotents. The 
automorphism group of the algebra often preserves the ``Frobenius strata'', i.e., the 
idempotents having the same Frobenius norm.
\begin{enumerate}
\item[(a)] Is there any natural description of the Frobenius strata for axial algebras?
\item[(b)] Do Frobenius strata have any distinguished algebraic structure?
\item[(c)] Are there any natural obstructions for fusion laws of fixed strata in (a)?
\end{enumerate}
\begin{flushright} V.G. Tkachev\end{flushright}
\end{prob}

Related to this question is the following. The automorphism group of $A$ acts on the space 
of all Frobenius forms defined on $A$. In most known examples, this space is $1$-dimensional 
and the forms are distinguished by the length of some specific invariant element of $A$, 
e.g., the identity element. Hence in the majority of examples the full automorphism group 
preserves the Frobenius form.

\begin{prob}
Does there exist an axial algebra with a Frobenius form such that the full automorphism group 
of the algebra does not preserve the form?
\begin{flushright} S. Shpectorov \end{flushright}
\end{prob}

Here we want a non-trivial example, excluding the likes of the following: take an algebra $B$ 
with a Frobenius form and define $A=B\oplus B$ with the Frobenius form on $A$ non-zero on the 
first summand and zero on the second one. Then this form will clearly not be preserved by the 
automorphisms switching the two summands $B$. 

Especially, we are interested in the case of algebras where the space of Frobenius forms is 
$1$-dimensional. Then every automorphism $g$ alters the Frobenius form by a scalar $\al_g$ and 
the map $g\to\al_g$ is a homomorphism from $G=\Aut(A)$ to $\F^\sharp$. For this homomorphism 
to be nontrivial, $A$ cannot be unital.

\subsection{Structure}

The concept of the radical of an axial algebra was introduced in \cite{kms}.

\begin{defn}
The \emph{radical} $R(A)$ of a primitive axial algebra $(A,X)$ is the unique largest
ideal of $A$ not containing any axes from $X$.
\end{defn}

For example, in a baric axial algebra, the baric radical is an ideal of codimension $1$ 
containing no axes from $X$, hence it coincides with the radical of the algebra.

When $A$ admits a Frobenius form, the radical $R(A)$ is always contained in the radical 
$A^\perp$ of the Frobenius form. Furthermore, $R(A)=A^\perp$ if and only if $(a,a)\neq 0$ 
for all generating axes $a\in X$. That is, all generating axes are non-singular for the 
Frobenius form.

The recent article \cite{msz} introduced the concept of the Jacobson radical of an axial 
algebra.

\begin{defn}
The \emph{Jacobson radical} $J(A)$ of an axial algebra $A$ is the intersection of all
maximal ideals of $A$. If $A$ has no maximal ideals then we simply take $J(A)=A$.
\end{defn}

It is shown in \cite{msz} that $R(A)\subseteq J(A)\subseteq A^\perp$. Clearly, when all axes 
$a\in X$ are non-singular, all three radicals coincide. However, as a general axial algebra 
$A$ may admit different Frobenius forms, including the ones not satisfying the non-singularity 
condition, we cannot expect that in general $R(A)$ and $J(A)$ coincide with $A^\perp$, but they 
can, in principle, be equal to each other.

\begin{prob}
Is it true that the axial radical $R(A)$ and the Jacobson radical $J(A)$ coincide in
every primitive axial algebra? Equivalently, is it true that $J(A)$ never contains
primitive axes?
\begin{flushright} A. Mamontov, S. Shpectorov, V. Zhelyabin \end{flushright}
\end{prob}

The paper \cite{msz} also introduces the following interesting concept.

\begin{defn} \label{block}
For an axial algebra $(A,X)$ and a primitive axis $a\in X$, the \emph{block} $I_a$ of $A$
corresponding to $a$ is defined as the ideal generated by $a$, i.e., the smallest
ideal of $A$ containing $a$.
\end{defn}

An algebra is said to be indecomposable when it is not the sum of its proper ideals.

\begin{prob}
Is it true that each axial block is indecomposable?
\begin{flushright} A. Mamontov, S. Shpectorov, V. Zhelyabin \end{flushright}
\end{prob}

On the set of blocks, we have the dominance relation defined by inclusion.

\begin{prob}
Can there be non-trivial embeddings between blocks, i.e., can the dominance relation
to be non-symmetric?
\begin{flushright} A. Mamontov, S. Shpectorov, V. Zhelyabin \end{flushright}
\end{prob}

No such examples are known at present. Furthermore, for some classes of axial 
algebras, e.g., for Majorana algebras, dominance is known to be symmetric.

\begin{defn}\label{non-ann}
Suppose that $(A,X)$ is a primitive axial algebra. The \emph{non-annihilation graph}
$\Dl$ of $A$ has $X$ as its set of vertices and distinct $a,b\in X$ are adjacent in the
non-annihilation graph if and only if $ab\neq 0$.  
\end{defn}

The following related question was first posed in \cite{kms} in the context of Majorana
algebras and then also repeated in a more general setting in the survey \cite{ms}.

\begin{prob}
Is it true that in the finest (direct) sum decomposition of a primitive axial algebra
the summands correspond to the connected components of the non-annihilation graph
$\Dl$?
\begin{flushright} J. McInroy, S. Shpectorov (Conjecture 3.16 [MS]) \end{flushright}
\end{prob}

The most interesting case of this is as follows.

\begin{prob}
Is it true that the non-annihilation graph $\Dl$ of a simple Majorana algebra is
necessarily connected?
\begin{flushright} A. Mamontov, S. Shpectorov, V. Zhelyabin \end{flushright}
\end{prob}

The great majority of the known examples of axial algebras are unital, that is, 
they contain a multiplicative identity. Furthermore, among the known examples, the algebras 
without identity are always non-simple. While it is unlikely that this is a general 
statement covering all primitive axial algebras, the following question is certainly 
an interesting one to investigate.

\begin{prob}
For which classes of (primitive) axial algebras are all simple algebras unital?
\begin{flushright} I. Gorshkov, S. Shpectorov \end{flushright}
\end{prob}

Also, the following is a related question.

\begin{prob} \label{trace}
Unital axial algebras $A$ with Frobenius form have a trace form, i.e., a distinguished linear map $Tr:A\to\F$ such that $(u,v)= Tr(uv)$ for all $u,v\in A$. To what extent is this true for non-unital axial algebras?
\begin{flushright} J. Desmet \end{flushright}
\end{prob}

\subsection{Inner Lie algebra}

The Lie algebra of inner derivations plays an important role in the theory of Jordan algebras. Can there be a similar theory for more general axial algebras?

Schafer \cite{s} defined the Lie algebra of inner maps of a non-associative algebra $A$ as the subalgebra of $\End(A)^\circ$ generated by all left and right adjoints in $A$. Here $\End(A)^\circ$ is the Lie algebra of all endomorphisms of the vector space $A$ with respect to the bracket operation. In our case, $A$ is commutative, so we do not distinguish between left and right adjoint maps. Furthermore, we will focus on the commutator subalgebra of Scharfer's Lie algebra of inner maps.

\begin{defn}\label{associator}
For a (commutative) non-associative algebra $A$ over a field $\F$, and $a, b \in A$, the \emph{associator} $D_{a,b}:A\to A$ is defined by $D_{a,b}(x)=a(bx)-b(ax)$. 
\end{defn}

Clearly, $D_{a,b}=\ad_a\ad_b-\ad_b\ad_a=[\ad_a,\ad_b]$.

When the algebra $A$ is Jordan, the associators $D_{a,b}$ are derivations of $A$ for all choices of $a,b\in A$. For algebras of Jordan type half, it was shown in \cite{d25} that $D_{a,b}$ are derivations for all pairs of axes $a$ and $b$, except when $a$ and $b$ generate the so-called \emph{Fischer triple}, a subalgebra isomorphic to $3C(\frac{1}{2})$ and containing exactly three axes, $a$, $b$ and $c=a^{\tau_b}=b^{\tau_a}$, even over any field extension. Furthermore, if $D_{a,b}$ are derivations for all pairs of axes $a,b$, then $A$ is Jordan (see \cite{d25}).

This discussion makes the following object the focus of study.

\begin{defn}\label{inner Lie algebra}
The algebra $DA_X=\lla D_{a,b}\mid a,b\in X\rra$, taken with the bracket multiplication, is called the \emph{inner Lie algebra} of the axial algebra $(A,X)$.
\end{defn}

Clearly, this is a subalgebra of the commutator subalgebra of Schafer's Lie algebra of inner maps on $A$.

\begin{prob} \label{g1}
Let $A$ be an axial algebra and $X$ and $X'$ be two generating sets of axes of $A$. Are the algebras $DA_X$ and $DA_{X'}$ equal? An important special case of this question arises when the fusion law of $A$ is graded and $X'=\bar X$ is the closure of the set $X$.
\begin{flushright} I. Gorshkov \end{flushright}
\end{prob}

In essence, the question asks whether the inner Lie algebra of an axial algebra is stable with respect to the closure of the set of axes, as defined in \cite{kms}, where it is shown that many other objects associated with an axial algebra are stable. 

\begin{prob}
For $(A,X)$ as in Problem $\ref{g1}$, if $G$ is an automorphism group of $A$ preserving $X$ then, clearly, $G$ acts as an automorphism group of $DA_X$. Is it true that the automorphism groups of the algebra $A$ and the algebra of inner mappings $DA_X$ for any generating set $X$ are isomorphic?
\begin{flushright} I. Gorshkov \end{flushright}
\end{prob}

When the set $X$ is invariant under the full automorphism group of $A$, e.g., if $X$ consists of all primitive axes found in $A$, we have that $\Aut(A)$ is contained in $\Aut(DA_X)$. 
The reverse inclusion follows if we can somehow recover $A$ from $DA_X$.

\begin{prob}   
Construct axial algebras such that the algebras of inner mappings are simple Lie algebras. E.g., for some classical Lie algebras and $G_2$, one can take a simple Jordan algebra. Can one construct such examples for the remaining simple Lie algebras, or can we construct examples with different fusion laws?
\begin{flushright} I. Gorshkov \end{flushright}
\end{prob}


More generally, it would be interesting to see how the structure of $A$ is reflected in its inner Lie algebra. In particular, the following question is interesting.

\begin{prob}   
Is it true that the simplicity of an axial algebra $A$ implies the simplicity of the inner Lie algebra $DA_X$ for any generating set of axes $X$? Vice versa, does 
simplicity of $DA_X$ imply simplicity of $A$?
\begin{flushright} I. Gorshkov \end{flushright}
\end{prob}

The notion of a conservative algebra is due to Kantor \cite{k72}.

\begin{defn}
An algebra $A$ is said to be \emph{conservative} if the brackets of two left adjoints is again a left adjoint, that is,
$$
[L_a,L_b]=L_{g(a,b)}
$$
for all $a,b\in A$ and some function $g:A\times A\to A$. Here we use the notation of $L_a$ for the left adjoint operator of $a\in A$.
\end{defn}

Associative and Lie algebras are conservative, and so are also Jordan algebras. In the context of axial algebras, we can write $[\ad_a,\ad_b]=\ad_{g(a,b)}$. The condition itself means that the space of (left) adjoints is closed under the bracket operation and hence is a Lie algebra.

\begin{prob}
Characterise conservative axial algebras.
\begin{flushright} I. Kaygorodov\end{flushright}
\end{prob}

\subsection{Gelfand-Kirillov dimension}

\begin{defn}
Let $A$ be a finitely generated algebra over $\F$. Let $V\subseteq A$ be a finite-dimensional subspace generating $A$ and let $V^n$ denote the linear span of all products of length $n$ in elements of $V$. The Gelfand-Kirillov dimension of $A$ is defined as:
$$
\mathrm{GKdim}(A) = \limsup_{n \to \infty} \frac{\ln(\dim V^n)}{\ln n}.
$$
\end{defn}

Note that $\mathrm{GKdim}(A)=0$ if and only if $A$ is finite-dimensional, hence it is only interesting for infinite-dimensional algebras. According to \cite{kl}, there are no algebras $A$ with $0<\mathrm{GKdim}(A)<1$. Furthermore, Bergman \cite{b} showed that there exists no associative algebra $A$ with $1<\mathrm{GKdim}(A)<2$ and for all $s\geq 2$ there exist associative algebras with $\mathrm{GKdim}(A)=s$. Martinez and Zelmanov \cite{mz} characterised Jordan algebras $A$ with $\mathrm{GKdim}(A)=1$ and also showed that there is no Jordan algebra $A$ with $1<\mathrm{GKdim}(A)<2$.

\begin{prob}
\begin{enumerate}
\item[(a)] Is there an axial algebra of Gelfand-Kirillov dimension $s$, $1 < s < 2$?
\item[(b)] Classify axial algebras of Gelfand-Kirillov dimension $1$.
\end{enumerate}
\begin{flushright} I. Kaygorodov\end{flushright}
\end{prob}

Currently, there are only a few examples of axial algebras of infinite dimension are known. Centrone \cite{ce} checked that the Gelfand-Kirillov dimension of the 
Highwater algebras $\cH$ (see \cite{fms}) is $1$. 

\begin{prob}
Find the Gelfand-Kirillov dimension of the cover $\hat\cH$ of the Highwater algebra $\cH$.
\begin{flushright} C. Franchi, M. Mainardis\end{flushright}
\end{prob}

The cover algebra $\hat\cH$ was constructed by Franchi, Mainardis and McInroy \cite{fmm}.

\begin{prob}
Find the Gelfand-Kirillov dimension of the baric near Jordan type algebras constructed by Shi \textup{\cite{ys}}. 
\begin{flushright} I. Gorshkov, S. Shpectorov\end{flushright}
\end{prob}

\section{Algebras of Jordan type} 

The class of algebras of Jordan type has been the focus of much attention.

\subsection{General questions}

First of all, let us examine the very definition of the class of algebras of Jordan type. Recall, that primitivity of the generating axes is postulated and this property is used heavily in the available proofs.

\begin{prob}
Can a theory of algebras with Jordan type fusion be developed without the assumption of primitivity?
\begin{flushright} H. Cuypers\end{flushright}
\end{prob}

One of the immediate victims of this change is that we cannot guarantee the existence of a Frobenius form on the algebra.

\medskip
Currently known examples of algebras of Jordan type $\eta$ include Matsuo algebras introduced in Definition \ref{Matsuo} in terms of $3$-transposition groups and all Jordan algebras generated by primitive idempotents. 

\begin{prob}   
Does there exist a purely `algebraic’ definition of a Matsuo algebra, not referring to a group structure? In other words, given a commutative non-associative algebra with no distinguished set of primitive axes, how can one verify whether this algebra is Matsuo?
\begin{flushright} V.G. Tkachev\end{flushright}
\end{prob}

The classical approach to non-associtive algebras is via polynomial identities.

\begin{prob}
Do all Matsuo algebras satisfy identities not implied by commutativity? If they do, what is the minimal degree of an  identity of all Matsuo algebras not implied by commutativity? 
\begin{flushright} L. Rowen\end{flushright}
\end{prob}

Rowen pointed out that for Matsuo algebras, which are not Jordan algebras, the degree would have to be at least $4$ by a result of Chayet and  Garibaldi \cite[Proposition A.8]{cg}. Note that this result refers to the paper by Osborn \cite{o} where unitality of the algebra is assumed. Matsuo algebras are generically unital, as for each $3$-transposition group $(G,D)$, there is only a finite number of values of $\eta$ for which the Matsuo algebra $M_\eta(G,D)$ is non-unital. It is natural to assume that whatever identity holds for the generic Matsuo algebras will also hold in the exceptional non-unital situations.

\medskip
Recall that, for $\eta\neq\frac{1}{2}$, every algebra of Jordan type $\eta$ is a Matsuo algebra or a quotient of a Matsuo algebra \cite{hrs2}. So $\eta=\frac{1}{2}$ is the case where the Jordan algebras appear and where any new examples may arise. Note that a general commutative non-associative algebra satisfying a single-variable 
homogeneous polynomial identity must have $\frac{1}{2}$ in its spectrum \cite{t21}.

\begin{prob}   
Do the axial algebras of Jordan type half satisfy a non-trivial single-variable identity? For example, every Jordan algebra is power associative, e.g., it satisfies $xx^3 = x^2x^2$. Does this property hold at least for the non-Matsuo algebras of Jordan type half?
\begin{flushright} V.G. Tkachev\end{flushright}
\end{prob}

Recall that all Matsuo algebras with $\eta=\frac{1}{2}$ that are Jordan were determined by De Medts and Rehren \cite{dr}. Gorshkov, Mamontov and Staroletov \cite{gms} also determined the Matsuo algebras having non-trivial Jordan quotients. It is easy to check that the smallest non-Jordan Matsuo algebra of Jordan type half 
(for the group $2^4:S_4$) is not power associative. So we can expect the broad majority 
of Matsuo algebras to be non-power-associative.

Segev \cite{se} investigated primitive axes of Jordan type half in commutative non-associative algebras $A$ satisfying the identities $x^2x^2=xx^3$ and $x^2x^3=xx^4$ strictly, that is, these identities hold in $A$ over an arbitrary field extension. In particular, he shows that power associative algebras of Jordan type half are Jordan.

\medskip
There is a number of interesting questions related to the geometry of primitive axes in algebras of Jordan type. When $\eta\neq\frac{1}{2}$ and the algebra is a factor of some Matsuo algebra $M_\eta(G,D)$, primitive axes always form a Fischer space of a potentially larger $3$-transposition group $(\hat G,\hat D)$. Naturally, in a great majority of cases, we have $\hat G=G$. However, Alharbi \cite{a} also found several series of examples, for specific $3$-transposition groups and specific values of $\eta$, where $\hat G>G$. This is related to the question of finding the full automorphism group of the algebra, as $\hat G$ will be a normal subgroup of the full automorphism group.

Alharbi's work \cite{a} only covers the irreducible case of $3$-transposition groups. In this respect the following questions is important.

\begin{prob}
Extend Alharbi's classification of Matsuo algebras containing additional primitive axes to cover all reducible $3$-transposition groups. Equivalently, find all situations where the full automorphism group of $M_\eta(G,D)$, $\eta\neq\frac{1}{2}$ does not coincide with the automorphism group of $(G,D)$ and its Fischer space.
\begin{flushright} S. Shpectorov \end{flushright}
\end{prob}

We believe that there could be a blanket result covering all groups of $3$-transpositions. It would also be very good to extend Alharbi's classification on all quotients of Matsuo algebras.

The case $\eta=\frac{1}{2}$ is of course also important, but it is significantly more complex, especially since we do not know all examples here. Here we may have analogues of Alharbi's examples, but there also are further examples, where the set of primitive axes becomes infinite. This happens, for example, in Jordan algebras, but there are also examples of non-Jordan Matsuo algebras (see e.g. Example 7.2 in \cite{gss}) where this is the case. This has to do with the concept of a solid subalgebra (a $2$-generated subalgebra where all primitive idempotents are axes), which we discuss later. It follows from the theory of solid subalgebras (see \cite{ds}) that every algebra of Jordan type half without infinite solid subalgebras is a Matsuo algebra or a quotient of a Matsuo algebra. Desmet \cite{d26} has determined all Matsuo algebras with infinite solid lines, so this gives us a classification of algebras of Jordan type half with a finite set of primitive axes.

\begin{prob}
Are there any examples of Matsuo algebra $M_\eta(G,D)$ for $\eta=\frac{1}{2}$, where the set of all primitive axes is finite, but we still have the full automorphism group larger than the automorphism group of $(G,D)$?
\begin{flushright} S. Shpectorov \end{flushright}
\end{prob}

Note that here the set of all primitive axes will also form the Fischer space of some larger $3$-transposition group $(\hat G,\hat D)$, so this can be investigated by the same methods employed by Alharbi.

To finish with the topic of automorphism groups of algebras of Jordan type, let us mention the following question.

\begin{prob} \label{groups of Jordan type}
Which finite groups arise as Miyamoto groups of algebras of Jordan type?
\begin{flushright} H. Cuypers \end{flushright}
\end{prob}

Clearly, every finite group can be embedded, say, into the symmetric group and made to act on a Jordan algebra. The interesting case is where the group is in fact generated by the Miyamoto involutions, and so it is the Miyamoto group of a finite closed set of axes. We can call a finite group a \emph{Jordan type group} if it arises this way. For example, $J_2$ is a group of Jordan type, as found computationally by McInroy and Shpectorov.

\medskip
Tkachev believes that it is interesting to investigate the geometry of all primitive indempotents (not just axes) in Matsuo algebras and more generally in algebras of Jordan type. The example to consider here is the $2$-generated algebras of Jordan type. 

Generically, if $A=\lla a,b\rra$ is a $2$-generated axial Jordan type $\eta$ algebra, then $\eta\neq\frac{1}{2}$, $A$ has dimension $3$ and it contains the identity element $\mathbbm{1}$. Then $A$ contains exactly eight idempotents:
\begin{itemize}
\item[(i)] $0$ with spectrum $0^3$;
\item[(ii)] $a$, $b$, and  $c := a^{\tau_b} = b^{\tau_a}$ with spectrum $0^1 1^1\eta^1$;
\item[(iii)] $\mathbbm{1}-a$, $\mathbbm{1}-b$, and $\mathbbm{1}-c$ with spectrum $0^11^1(1 - \eta)^1$; and 
\item[(iv)] $\mathbbm{1}$ with spectrum $1^3$.
\end{itemize}
The idempotents in (ii) and (iii) are the primitive ones and they are all axes for different but similar fusion laws. Can a classification of this sort be developed for larger algebras?

\begin{prob}   
\begin{enumerate}
\item[(a)] Describe or characterise in some compact terms the set of all idempotents in a generic axial algebra of Jordan type.
\item[(b)] If $A$ is $m$-generated, and $C(m)$ is the maximal possible dimension, what can be said in the case where $\dim A = C(m)$.
\end{enumerate} 
\begin{flushright} V.G. Tkachev\end{flushright}
\end{prob}

Part (b) of this question assumes the existence of a function $C(n)$ bounding the dimension of an algebra of Jordan type with 
a given number $n$ of generators. This brings us to one of the central topics in the theory of algebras of Jordan type. 

\subsection{Small number of generators}

The existence of the function $C(n)$ bounding the dimension of the $n$-generated algebra of Jordan type is unknown for 
a general $n$, and this is the point of the following question.
 
\begin{prob} \label{Cm}
Is it true that every finitely generated algebra of Jordan type has finite dimension? If so, is the dimension bounded by a function of the number of generators?
\begin{flushright} Y. Segev\end{flushright}
\end{prob}

Finiteness of dimension is known for small values of $n$. For $n=2$, it was shown in \cite{hrs2}, where algebras of Jordan type were introduced, that every 
such $2$-generated algebra has dimension at most $3$. (That is, $C(2)=3$.) In the key case of $\eta=\frac{1}{2}$, there exists a $1$-parameter family $S(\frac{1}{2},\phi)$ 
of $3$-dimensional examples and every $2$-generated algebra of Jordan type half is a quotient of one of these examples. Similarly, for $n=3$, it was shown in \cite{gs} 
that $C(3)=9$. All $9$-dimensional $3$-generated algebras of Jordan type half form a $4$-parameter family, and every $3$-generated example is a quotient of one of the 
$9$-dimensional ``universal'' algebras. For $n=4$, it was shown in \cite{drs} that $C(4)=81$ and that this bound is exact. However, there is no known ``universality'' 
statement about $4$-generated algebras of Jordan type. So far as we know, no results are available for $n\geq 5$.

In view of this, it would be interesting to classify algebras of Jordan type with a small number of generating axes satisfying additional restrictions. 
Define the \emph{non-orthogonality} graph of the set $X$ of generating axes as follows. The vertices of the graph are the axes from $X$ and two axes 
$a,b\in X$ are adjacent when $(a,b)\neq 0$. This concept first appeared in the structure theory of axial algebras.

\begin{prob} 
Classify algebras $A$ of Jordan type $\frac{1}{2}$ whose non-orthogonality graph $\Gm$ is: (a) a star; (b) a chain; (c) a tree?
\begin{flushright} I. Gorshkov \end{flushright}
\end{prob}

Note that the non-annihilation condition $ab\neq 0$ is a weaker form of the condition $(a,b)\neq 0$. It would also be interesting to classify 
algebras of Jordan type with a given non-annihilating graph (see Definition \ref{non-ann}) on the set of generators.

Here it is important that the set $X$ of generating axes is not closed and, in some sense, minimal.

Let us also formulate a few related questions.

\begin{prob}
Let $A$ be a finitely generated algebra of Jordan type half. What is the maximum dimension of the radical of $A$?
\begin{flushright} I. Gorshkov\end{flushright}
\end{prob}

Of course, we have examples of baric algebras of Jordan type, where the radical is of codimension $1$. Do these baric examples provide the upper 
bound on the dimension of the radical?

The following property is also related to the question of finiteness of the algebra dimension.

\begin{defn}
For an arbitrary algebra $A$ and a generating set $X$ of $A$, the \emph{length} of $X$ is defined as the minimum $k$ 
such that every element of $A$ is a linear combination of products of elements from $X$ with at most $k$ factors.

The \emph{length} of $A$ is the maximum length of $X$ over all generating sets $X$ of $A$.
\end{defn}

Clearly, if the algebra is finitely generated and the length is finite then the algebra is finite-dimensional.

\medskip
Length of non-associative algebras was studied by Guterman and co-authors (see, e.g., Guterman and Kudryavtsev \cite{gk}).
In the context of axial algebras, it is natural to restrict generating sets $X$ to just the generating sets of 
axes.

\begin{prob}
Give an upper bound on the (axial) length of finite-dimensional (primitive) axial algebras.
Calculate the length of axial algebras of Jordan type.
\begin{flushright} I. Kaygorodov\end{flushright}
\end{prob}

Let $L_k$ be the subspace of all elements of $A$ that are linear combinations of products of elements of $X$ with at most $k$ factors. Then the sequence $L_0\subseteq L_1\subseteq L_2\subseteq\ldots\subseteq L_i\subseteq\ldots$ is a filtration of $A$, and hence we can define a graded algebra structure on
$$
\hat A=L_0\oplus L_1/L_0\oplus L_2/L_i\oplus\ldots L_i/L_{i-1}\oplus\ldots
$$

\begin{prob}
Is the associated graded algebra of a finitely-generated algebra of Jordan type half Jordan? Is it finite-dimensional? 
\begin{flushright} L. Rowen\end{flushright}
\end{prob}

Note that the dimension of the graded algebra $\hat A$ has the same dimension as $A$. If $\hat A$ happens to be Jordan and if this means that $\hat A$ is finite-dimensional for a finite generating set $X$, this also implies that finitely generated algebras of Jordan type half are finite-dimensional. Note that Rowen suggested that the answers to the above questions may be different in small field characteristic, such as $3$ and $5$. 

Following \cite{hrs1}, one can define the universal $k$-generated algebra of Jordan type half endowed with a Frobenius form such that $(a,a)=1$ for each generating axis $a$.
This universal algebra is defined over a ring, whose spectrum represents the variety of all $n$-generated algebras of Jordan type half over fields. For $n=2$, this variety is an affine line representing the $1$-parameter family of $3$-dimensional examples plus a few isolated points for the occasional $2$-dimensional factors. For $n=3$, the variety includes
a $4$-dimensional affine space for the $4$-parameters family of $9$-dimensional examples plus additional components for the available quotients.

\begin{prob}
Describe the complete structure of the variety of $3$-generated algebras of Jordan type.
\begin{flushright} I. Gorshkov and S. Shpectorov\end{flushright}
\end{prob}

This involves a complete analysis of all values of the four parameters, for which the $9$-dimensional algebra is not simple and hence has quotients. A large amount of work in this direction was done by Bildanov and Gorshkov \cite{bg}, where they analyse 
simple factors of $9$-dimensional universal algebra, assuming that the field is quadratically closed.

\medskip
A much harder project would be to analyse to a similar level of detail the case of $n=4$.

\begin{prob}
Analyse the variety of algebras of Jordan type for $n=4$. Is there a single multi-dimensional component in this variety, from which all or almost 
all other algebras can be obtained as quotients. If such a principal component exists then what is its dimension (number of parameters describing the algebras) and what is the dimension 
of the algebras in this family?
\begin{flushright} I. Gorshkov and S. Shpectorov\end{flushright}
\end{prob}

Shi \cite[Chapter 9]{ys} analysed the case of $n=4$, where we additionally assume that $(a,b)=1$ for all $a,b\in X=\bar X$ (equivalently, the algebra 
is baric). What she found is that the universal algebra is $81$-dimensional only when the field is in characteristic $3$ and in all other characteristics 
the dimension of the universal algebra is $54$. We can hence conjecture that the dimension of algebras in the principal component is $54$. 

All algebras of Jordan type for $n=2$ and $3$ are Jordan. For $n=4$, we get the first examples of non-Jordan algebras of Jordan type. It seems possible 
that the bulk of all examples for $n=4$ or any larger $n$ consists of the Jordan algebras described by a single principal component, as in the above 
problem, while the non-Jordan examples provide additional components of smaller component dimension.

Related to the above discussion is the following question.

\begin{prob}
What is the radical of the universal $n$-generated algebra of Jordan type half?
\begin{flushright} L. Rowen\end{flushright}
\end{prob}

Note that the universal algebra has lots of quotients and hence also lots of ideals not containing axes. However, as it is defined over a ring, the concept of radical needs to be examined to see if it still makes sense.

Rowen \cite{ro} also developed a theory of generic (universal) algebras of Jordan type. This paper goes a long way towards clarifying the structure of the variety of algebras of Jordan type. In particular, he constructed several universal objects, one representing the Matsuo branch, one the Jordan branch and several more that are neither Matsuo, nor Jordan.

\medskip
We wrap this subsection with two further questions from Rowen.

\begin{prob}
In an algebra of Jordan type half, does $(x,y) = 0$ imply $xy$ is $4$-nilpotent, i.e., $(xy)^2 = 0$ or $(xy)(xy)^3 = 0$?  (These make sense even when power nilpotence fails.) More generally, which identities can be deduced from orthogonality?
\begin{flushright} L. Rowen\end{flushright}

\end{prob}
\begin{prob}
Can there be a non-computer proof that $3$-generated algebras of Jordan type are Jordan? 
\begin{flushright} L. Rowen\end{flushright}
\end{prob}

The available computer-aided proof in \cite{gs} uses $GAP$ to verify the expressions for specific products and also to check that the resulting $9$-dimensional algebra is always Jordan. So the computer is only used for routine calculations, which can in principle also be done by hand. So the above question can perhaps be interpreted as follows: can the Jordan property be shown directly, without first reconstructing the whole algebra? We note that Rowen \cite{ro} also provides an answer to this question when the characteristic of the field is different from $3$ and $5$.

\subsection{Modules for algebras of Jordan type}

Let $A$ be an algebra of Jordan type $\eta$ over a field $\F$.

\begin{defn}
An $A$-module is a vector space $M$ over $\F$ together with 
a product $A\times M\to M$ such that: 
\begin{itemize}
\item[(a)] $M=M_1(a)\oplus M_0(a)\oplus M_\eta(a)$; 
and 
\item[(b)] the Jordan type  fusion law, 
$A_\lm(a)M_\mu(a)\subseteq M_{\lm\ast\mu}(a)$, is satisfied for each 
axis $a\in A$. 
\end{itemize}
The module $M$ is called \emph{pure} if $M_1(a)=0$.
\end{defn}

Here, as usual, $M_\lm(a)=\{u\in M\mid au=\lm u\}$ and 
$M_\Lm(a)=\oplus_{\lm\in\Lm}M_\lm(a)$.

Note that this is a particular case of a more general concept of a module for an axial algebra as introduced by De Medts and Van Couwenberghe \cite{dv}, although our terminology is slightly different.

\begin{exa}
If $A$ is included in a larger algebra $B$ of the same Jordan type then 
$A^\perp=\{u\in B\mid(u,A)=0\}$ is a pure $A$-module.
\end{exa}

\begin{prob}
Develop the theory of pure $A$-modules. Can they be classified, say, 
for the $2$-generated algebras of Jordan type?
\begin{flushright} S. Shpectorov \end{flushright}
\end{prob}

Because of the above example, modules can be used in the structure theory of algebras of Jordan type. Clearly, a direct sum of modules 
is again a module, and so we need to focus on indecomposable modules.

\begin{prob}
Is it true that the dimension of an indecomposable module $M$ is bounded by a function of the number of generating axes of the algebra $A$?
\begin{flushright} S. Shpectorov \end{flushright}
\end{prob}

This question is motivated by Segev's Question \ref{Cm}.

\subsection{Solid subalgebras, components}

We have already mentioned solid subalgebras. Here is the exact definition. 

\begin{defn}
Suppose that $(A,X)$ is an algebra of Jordan type half and $a,b\in X$. The subalgebra $B=\lla a,b\rra$ is said to be \emph{solid} if every primitive idempotent in $B$ is a primitive axis in the entire algebra $A$. Furthermore, this must remains true over every extension of the ground field $\F$. 
\end{defn}

It was shown in \cite{gss,d25} that a $2$-generated subalgebra $B$ is solid unless $B\cong 3C(\frac{1}{2})$ and exactly three idempotents from $B$, namely, $a$, $b$ and $c=a^{\tau_a}=b^{\tau_a}$, are primitive axes of $A$. This exceptional configuration with just three axes is called the \emph{Fischer triple}. Note that, while the Fischer triple is isomorphic to $3C(\frac{1}{2})$, it does not mean that every such subalgebra is a Fischer triple: subalgebras $3C(\frac{1}{2})$ can also be solid. The difference appears when we extend the ground field. The Fischer triple has still just the three primitive axes $a$, $b$, and $c$, while the solid $3C(\frac{1}{2})$ contains more and more primitive idempotents and hence primitive axes, as we extend the field.

If the algebra $A$ contains no solid $2$-generated subalgebras other than $2B$ then every $2$-generated subalgebra in $A$ is either $2B$ or it is a Fischer triple $3C(\frac{1}{2})$. This implies that $A$ is a Matsuo algebra or a quotient of a Matsuo algebra. Because of this, we can focus on the algebras $A$ containing solid $2$-generated subalgebras that are not $2B$. Furthermore, we can assume that $\F$ is big enough, so that all such solid subalgebras contain at least four axes and can therefore be distinguished from the subalgebras $2B$ and Fischer triples $3C(\frac{1}{2})$.

Desmet \cite{d25} provided an alternative characterisation of solid subalgebras. Namely, he proved that $\lla a,b\rra$ is solid if and only if the associator $D_{a,b}$ (see Definition \ref{associator}) is a derivation of the algebra $A$. Note that for $\lla a,b\rra\cong 2B$ this condition is degenerate, as in this case $D_{a,b}=0$. When $A$ contains solid lines other than $2B$, at least some inner maps $D_{a,b}$ are non-zero derivations. We will call such $D_{a,b}$ \emph{inner derivations}. Naturally, the inner derivations of $A$ generate a subalgebra in the Lie algebra of all derivations of $A$.

\begin{prob}
Describe the cases where the algebra of inner derivations is not the entire Lie algebra of derivations. Of particular interest are the Matsuo algebras for which this holds.
\begin{flushright} J. Desmet, S. Shpectorov \end{flushright}
\end{prob}

In particular, this seems to happen over fields of characteristic $3$, when even in the $2$-generated case, the Lie algebra of all derivations is bigger than the $1$-dimensional subalgebra of inner derivations.

Recall the concept of the inner Lie algebra $DA_X$ of an axial algebra $A$ from Definition \ref{inner Lie algebra}.

\begin{prob}
What properties does the Lie algebra $DA_X$ have in the case where $A$ is an algebra of Jordan type $\frac{1}{2}$?
\begin{flushright} I. Gorshkov \end{flushright}
\end{prob}

It is tempting to ask if the definition of solid subalgebras can be generalised beyond the $2$-generated case. 

\begin{prob}
Suppose $B$ is a subalgebra of an algebra $A$ of Jordan type half and each primitive idempotent in $B$ is a 
primitive axis in the entire $A$. What can we say about $B$?
\begin{flushright} I. Gorshkov \end{flushright}
\end{prob}

As we will see shortly, such a subalgebra $B$ is necessarily Jordan, but can we say more? In particular, does this provide any structural information about $A$?

\medskip
Let us now discuss the concept of a component in an algebra $A$ of Jordan type half, stemming from \cite{ds}. All $2$-generated subalgebras in $A$ can be classified into several types:
\begin{enumerate}
\item[(a)] \emph{toric} subalgebras $3C(\phi)$, $\phi\neq 1,0$; these are solid, so we exclude the Fischer triples arising for $\phi=\frac{1}{2}$;
\item[(b)] Fischer triples $3C(\frac{1}{2})$;
\item[(c)] \emph{baric} subalgebra $3C(1)$;
\item[(d)] its $2$-dimensional \emph{flat} baric quotient $\overline{3C(1)}$;
\item[(e)] \emph{baric pair} $C(0)$; and finally,
\item[(f)] its $2$-dimensional quotient $2B$, which we will occasionally refer to as an \emph{orthogonal pair}.
\end{enumerate}

Recall that the term baric above refers to the property of the axial algebra to have a radical $R(A)$ of codimension $1$, or equivalently, having an algebra homomorphism $A\to\F$ mapping all primitive axes to $1$. The algebra $3C(0)$ is not itself baric, but it is made out of two flat baric subalgebras intersecting in a $1$-dimensional ideal, the radical of $3C(0)$.

\begin{defn}
Given a solid subalgebra $B=\lla a,b\rra$ in an algebra $A$ of Jordan type half, we say that $B$ is \emph{connected} if and only if $(a,b)\neq 0$.
\end{defn}

This refers to connectivity (or rather irreducibility) of the variety of primitive idempotents in $B$. 
Note that the only disconnected subalgebras $B$ are those where $(a,b)=0$ and where $(a,b)=\frac{1}{4}$ and $B$ is a Fischer triple.

In order to introduce components, let us first discuss the generating set of axes $X$ of $A$. We will assume that $X$ is closed, that is, $X$ is an axet. However, this is not enough as we want to include the additional primitive axes arising in solid lines. 

\begin{defn}
An axet $X$ (i.e., a closed generating set of axes) in an algebra $A$ of Jordan type half is said to be \emph{solid} if, for all $a,b\in X$, the intersection $B\cap X$ is either a Fischer triple or it includes all primitive axes from $B=\lla a,b\rra$. 
\end{defn}

That is, a solid axet is also closed with respect to solid subalgebras. As before, we may assume that the ground field $\F$ is big enough so that the solid subalgebras $B$ can already be distinguished from the Fischer triples. 

From now on we assume that $X$ is a solid axet in $A$.

\begin{defn}
The \emph{connectivity graph} of $(A,X)$ has $X$ as its set of vertices. Two axes $a,b\in X$ are connected by an edge if $B=\lla a,b\rra$ is connected. 
\end{defn}

Note that \cite{ds} has a slightly different definition, where $a$ and $b$ are connected by an edge if $a$ and $b$ are contained in a connected solid $2$-generated subalgebra $B$. This definition is slightly better for proofs, but it leads to exactly the same concepts of components.

\begin{defn}
We call the connected components of the connectivity graph the \emph{components} of $X$. Furthermore, the subalgebras of $A$ generated by the components of $X$ will be called the \emph{component subalgebras} of $A$.
\end{defn}

It is shown in \cite{ds} that each component subalgebra of an algebra $A$ of Jordan type half is a Jordan algebra, so they may also be referred to as \emph{Jordan components} of $A$. It is also shown in \cite{ds} that the set of all Jordan components admits a natural structure of a Fischer space, and so there is a $3$-transposition group associated with each $A$. 
Hence, Matsuo and Jordan algebras are simply two extreme cases of algebras of Jordan type: Jordan algebras arise when there is only one component of an arbitrary size, and Matsuo algebras appear when we have an arbitrary Fischer space of components, but they are all $1$-dimensional.

Of course, the most interesting and least explored case is the one between the two extremes.

\begin{prob}
Determine which Jordan algebras can arise as the Jordan components of an algebra of Jordan type $\frac{1}{2}$. Does the structure of the Jordan components depend on the Fischer space?
\begin{flushright} J. Desmet, S. Shpectorov \end{flushright}
\end{prob}

Clearly, every Jordan algebra generated by primitive idempotents can be the whole algebra $A$. So our question is primarily about the situation where we have more than one component, that is, the Fischer space of components is non-trivial.

Suppose that $J$ and $K$ are two Jordan components of the algebra $A$. It can be shown that the value $(x,y)\in\{0,\frac{1}{4}\}$ is independent of the choice of primitive axes $x\in J$ and $y\in K$. Hence each component $J$ has a hyperplane that is orthogonal to each component $K\neq J$. In this context, the next question is very interesting. 

\begin{prob}
Is it true that the radical of a component is in the radical of the entire algebra?
\begin{flushright} S. Shpectorov \end{flushright}
\end{prob}


Experimenting with small component algebras, Desmet constructed examples of algebras of Jordan type half that are neither Jordan, nor Matsuo, thus giving a negative answer to the original conjecture that every algebra of Jordan type is either Matsuo (or a factor of Matsuo algebra) or Jordan. In his examples, the components are baric, that is, they have a radical of codimension $1$. Furthermore, the sum of the radicals of all components comprise the radical of $A$ and modulo the radical the algebra becomes a Matsuo algebra.

In this respect, let us mention a related question. A \emph{central extension} of an algebra $A$ is an algebra $\hat A$ with an annihilator 
ideal $I$ (i.e., $AI=0$) such that $\hat A/I\cong A$.

\begin{prob}
Is it true that every algebra of Jordan type $\frac{1}{2}$ that is a central extension of a Jordan algebra is itself a Jordan algebra?
\begin{flushright} I. Kaygorodov, C. Martín González, P. Páez-Guillán \cite{kmp}\end{flushright}
\end{prob}

\subsection{Unitality}

It has been observed that every simple Matsuo algebra is unital, that is, it contains an identity. Also, every simple Jordan algebra, 
except the trivial $1$-dimensional one, contains an identity.

\begin{prob} \label{unital}
Is it true that every finitely generated simple axial algebra of Jordan type $\frac{1}{2}$ is unital?
\begin{flushright} I. Gorshkov \end{flushright}
\end{prob}

This is indirectly related to the following question, a special case of Problem \ref{trace}, which also contains the definition of a trace.

\begin{prob}
Does every algebra of Jordan type have a trace map?
\begin{flushright} J. Desmet \end{flushright}
\end{prob}

Indeed, if Problem \ref{unital} is answered in affirmative, then a trace map can be constructed on every algebra $A$ of Jordan type half. Indeed, consider a maximal ideal $I$ of $A$. Then $A/I$ is simple and hence unital by Problem \ref{unital}. Hence $A/I$ admits a trace map and it clearly lifts to the entire $A$. 


This construction also shows that the trace map does not have to be unique, especially when $A$ has a number of different simple factor algebras.

\subsection{Simplicity of the Miyamoto group}

If $A$ is a finite-dimensional simple Jordan algebra over a field of characteristic not $2$ then the Miyamoto group of $A$ is typically either simple or nearly simple.

\begin{prob}
When does the Miyamoto group of an algebra $A$ of Jordan type $\frac{1}{2}$ is simple or nearly simple? When this happens, is it true that $A$, as a module for its Miyamoto group, is the sum of the $1$-dimension module generated by the identity element and an irreducible complement?
\begin{flushright} I. Gorshkov \end{flushright}
\end{prob}

It is also interesting to ask the same questions about the subalgebra of an algebra of Jordan type half arising on the $0$-eigenspace $A_0(a)$ with respect to a primitive axis $a$.

\begin{prob}
Let $a$ be a primitive axis of an algebra $A$ of Jordan type $\frac{1}{2}$. When is $A_0(a)$ an axial algebra? When is it a simple axial algebra?
\begin{flushright} I. Gorshkov \end{flushright}
\end{prob}

\section{Algebras of Monster type} 

In this section we assembled questions about the larger and less investigated class of axial algebras, algebras of Monster type.

\subsection{General algebras of Monster type}

We had a fair number of questions earlier concerning the algebra of inner mappings of general axial algebras. 
Of course, these questions are especially interesting for important special classes of axial algebras like the algebras of Monster type.

\begin{prob}
Describe the Lie algebra of inner mappings of an algebra of Monster type.
\begin{flushright} I. Gorshkov \end{flushright}
\end{prob}

Naturally, it is interesting when the algebra of inner mappings consists of derivations.

\begin{prob}
In which algebras of Monster type $(\al,\bt)$ are all internal mappings $D_{a,b}$, with generating axes $a,b$, derivations? 
\begin{flushright} I. Gorshkov \end{flushright}
\end{prob}

By a result from \cite{a} and \cite{gmms} (see the proof of Theorem 3.1), $\frac{1}{2}$ must be an eigenvalue, so 
$\frac{1}{2}\in\{\al,\bt\}$ is a necessary condition. It is interesting to investigate whether there is any significant 
difference between the cases where $\al=\frac{1}{2}$ and where $\bt=\frac{1}{2}$. Based on the current, very scarce 
set of available examples, the following question deserves attention.

\begin{prob}
Are there any algebras of Monster type $(\al,\bt)$ with $\beta=\frac{1}{2}$ that admit non-trivial derivations? 
\begin{flushright} S. Shpectorov\end{flushright}
\end{prob}

We mean here the true examples of algebras of Monster type, not the algebras of Jordan type, where 
the eigenvalue $\eta=\frac{1}{2}$ can be treated as either $\al$ or $\bt$.

\medskip
All known examples of algebras of Monster type admit a Frobenius form. 

\begin{prob}
Is there always a Frobenius form on an algebra of Monster type?
\begin{flushright} S. Shpectorov \end{flushright}
\end{prob}

Finally, let us mention some more technical questions that arose from attempts to investigate 
the known examples of algebras of Monster type.

Recall that the fusion law of Monster type is Seress, that is, $0\ast\lm\subseteq\lm$ for all 
$\lm$. One consequence of this, known as the Seress Lemma, is that every axis $a$ of Monster 
type, whether primitive or not, associates with the corresponding eigenspace $A_0(a)$. Before Seress' 
observation, this property was deduced for the class of Majorana algebras from the Norton inequality:
$$
(uv,uv)\leq(u,u)(v,u),
$$
for all $u,v\in A$. One consequence of the Seress property is that $A_0(a)$ is a subalgebra, and this 
seems to hold more generally in algebras of Monster type. Namely, within the automorphism groups 
project \cite{gmms}, it was observed that, for all idempotents $b$ encountered in the calculations, 
the eigenspace $A_0(b)$ is always a subalgebra, regardless of the fusion law $b$ satisfies or even regardless of whether the adjoint of $b$ is semi-simple. Note 
that the algebras studied within this project are mostly Majorana algebras and they satisfy the Norton 
inequality. This leads to the following question.

\begin{prob}
Is it true that in an algebra $A$ with a Frobenius form satisfying the Norton inequality the eigenspace 
$A_0(b)$ is a subalgebra for each idempotent $b\in A$? Is this true in Majorana algebras or, more 
generally, in algebras of Monster type?
\begin{flushright} T.M. Mudziiri Shumba, S. Shpectorov\end{flushright}
\end{prob}

\subsection{$2$-Generated algebras and the $(\al,\bt)$-map}

The classification of $2$-generated algebras of Monster type has been the central research topic for the field of axial algebras. It started from the classification by Sakuma \cite{s} of OZ vertex operator algebras generated by two Ising vectors. This lead Ivanov \cite{i} to define the class of Majorana algebras and reprove Sakuma's result in this more general context in \cite{ipss}. Hall, Rehren and Shpectorov \cite{hrs1} generalised this result further 
to arbitrary algebras of Monster type $(\frac{1}{4},\frac{1}{32})$ 
admitting a normalised Frobenius form. Rehren \cite{r,r2} was the first one attempting to classify $2$-generated algebras of Monster type for arbitrary $(\al,\bt)$. In particular, he defined the universal $2$-generated algebra of Monster type without the Frobenius form assumption and used it to show that a $2$-generated algebra of Monster type is of dimension at most $8$, except possibly when $\al=2\bt$ or $\al=4\bt$. He also generalised the known Norton-Sakuma algebras to arbitrary $(\al,\bt)$. Later Franchi, Mainardis and Shpectorov simplified and slightly corrected Rehren's arguments in \cite{fms3} setting up a general classification of 
$2$-generated algebras. In this paper they also obtained the ultimate result about the Monster case $(\frac{1}{4},\frac{1}{32})$ showing that even after dropping all additional assumption the only examples are the Norton-Sakuma algebras. Furthermore, they classified the generic $2$-generated algebras of Monster type, that is, the ones that exist for all values of $(\al,\bt)$. In three spin-off articles \cite{fms,fms2,fm2} they described the only infinite-dimensional $2$-generated algebra of Monster type, called the Highwater algebra, which arises in one of the exceptional case of Rehren's, where $(\al,\bt)=(2,\frac{1}{2})$. They also completed the other exceptional case of type $(2\bt,\bt)$. At about the same time as these papers were written, Yabe \cite{y23} used new ideas to 
achieve and almost complete classification of symmetric $2$-generated algebras of arbitrary type $(\al,\bt)$. A $2$-generated algebra $\lla a,b\rra$ is \emph{symmetric} if it admits an automorphism swapping the generating axes $a$ and $b$. He also independently constructed the Highwater algebra and found several new families of $2$-generated algebras of Monster type. Franchi, Mainardis and McInroy \cite{fm2} completed the remaining open case of characteristic $5$ of symmetric $2$-generated algebras of Monster type and also classified factors of the Highwater algebra, which is baric and has a great many factors. Turner \cite{t,t2} started the classification of non-symmetric algebras of Monster type, and finally, the recently appeared monumental paper \cite{fmmt} completed the non-symmetric case and the entire classification of $2$-generated algebras of Monster type. 

In view of this major success, the focus now shifts from the classification of $2$-generated algebras to studying their properties and exploiting the classification. The questions about this can be formulated in geographical terms. Rehren \cite{r} was the first one to present the map of the $(\al,\bt)$ plane, where he indicated the points (types) for which various of his generalised Norton-Sakuma algebras exist, the type $(\frac{1}{4},\frac{1}{32})$ being the only point in the $(\al,\bt)$ plane for which all of these algebras are defined. After Sakuma found his families of $2$-generated algebras the map became a lot more complex, and the type $(\frac{1}{4},\frac{1}{32})$ is no longer the onely most interesting point in the plane.

\begin{prob}
Investigate the geography of the $(\al,\bt)$ plane and explicitly determine all types $(\al,\bt)$ for which many different $2$-generated algebra of Monster type exist.
\begin{flushright} S. Shpectorov \end{flushright}
\end{prob}

It is also interesting to determine points $(\al,\bt)$ for which algebras with specific properties may occur. 

\begin{prob}
Is it true that baric algebras of Monster type can only exist for $(\al,\bt)=(2,\frac{1}{2})$?
\begin{flushright} S. Shpectorov \end{flushright}
\end{prob}

Clearly, if $A$ is a baric axial algebra then all its axial subalgebras are also baric. In particular, all $2$-generated subalgebras in $A$ must be baric. We are only 
aware of the baric Highwater algebra among the $2$-generated algebras of Monster type. Hence this question requires checking that this is the only baric example on the list. 
It is not inconceivable that some exceptional baric algebras may exist in some positive characteristics.

\medskip
It has been noted that many interesting points in the $(\al,\bt)$ plane (the Monster type $(\frac{1}{4},\frac{1}{32})$, the Highwater type $(2,\frac{1}{2})$, and others) involve values that are powers of $2$.

\begin{prob}
Is there is systematic reason why $\al$ and $\bt$ that are powers of $2$ make more interesting types? 
\begin{flushright} C. Franchi, M. Mainardis \end{flushright}
\end{prob}

For example, the double axis construction for Jordan type suggests a possible reason, although it does not explain higher powers of $2$.

\medskip
Finally, we can look at the $(\al,\bt)$ plane from the group theory point of view.

\begin{prob}
Which finite groups can arise as Miyamoto groups of algebras of Monster type?
\begin{flushright} H. Cuypers \end{flushright}
\end{prob}

This is, of course, much harder than the corresponding question (Problem \ref{groups of Jordan type}) about algebras of Jordan type. This also needs to be studied according to the topology of the $(\al,\bt)$ plane, because some examples can arise for specific points $(\al,\bt)$, while others may fall into $1$- or $2$-parameter families.

\subsection{Larger algebras of Monster type}

Let $\F$ be the function field $\mathbb{Q}(\al,\bt)$. 

\begin{defn}
Algebras of Monster type $(\al,\bt)$ over $\F$ are called \emph{generic} algebras of Monster type.
\end{defn} 

The idea behind this definition is that by specialising $\al$ and $\bt$ we can construct a version of generic algebra for each choice of concrete values $\al$ and $\bt$. So these 
are the algebras that exist for all choices of the type $(\al,\bt)$. The Matsuo algebras and algebras $3A(\al,\bt)$ are the only currently known generic algebras of Monster type. More in particular, the only generic $2$-generated algebras are $2B$, $2A=3C(\al)$, $3A=3A(\al,\bt)$ and $3C=3C(\bt)$. In particular, the Miyamoto groups for generic algebras of Monster type are groups of $3$-transpositions.

\begin{prob}
Is it possible to classify all generic algebras of Monster type? In particular, are there any further generic algebras (other than the 
Matsuo algebras) for groups $S_n$, $n\geq 4$?
\begin{flushright} S. Shpectorov \end{flushright}
\end{prob} 

The case of the symmetric group is the natural first case to consider. In this case, the shape $2B3C$ clearly leads to Matsuo algebras, so the question is whether we can find 
any generic algebras for one of the other three shapes, $2A3A$, $2A3C$, and $2B3A$. It is known that the shape $2A3C$ for $(\al,\bt)=(\frac{1}{4},\frac{1}{32})$ collapses starting 
from the group $S_7$, so here we cannot have a full series of generic algebras. Also, it follows from Lemma 3.2 in \cite{fm} that the shape $2B3A$ collapses starting from $n=9$, at least in the class of Majorana algebras.

\begin{prob}
Find an upper bound on the dimension of a $3$-generated algebra of Monster type.
\begin{flushright} I. Gorshkov \end{flushright}
\end{prob}

The question is especially interesting for $(\al,\bt)=(\frac{1}{4},\frac{1}{32})$. The natural conjecture is that in this case the Griess algebra has the largest possible dimension, excluding Matsuo algebras and direct sums. In view of the infinite-dimensional Highwater algebra example, we need to avoid $(\al,\bt)=(2,\frac{1}{2})$.

\begin{defn}
An axial algebra $(A,X)$ is \emph{$n$-closed} if it is linearly spanned by the products of axes from $X$ with at most $n$ factors. 
\end{defn}

For example, it is known from \cite{hss1} that every algebra of Jordan type is $1$-closed as long as $X$ is a closet set of generating axes.

\begin{prob}
Is there a number $n$ such that every algebra of Monster type is at most $n$-closed?
\begin{flushright} I. Gorshkov \end{flushright}
\end{prob}

Naturally, we should assume that the set of axes is closed. Otherwise, we may need to deal with the following example: the Griess algebra is $3$-generated and for the set $X$ of size $3$ the value of $n$ will be quite large, at least $11$ and likely more. However, if we take a closed set then the Griess algebra is $1$-closed. 

Currently, we have examples of algebras that are $5$-closed with respect to a closed set of axes (axet). However, in all these algebras there is a second orbit of axes that brings the value of $n$ down.

\begin{prob}
Classify $1$-closed ($2$-closed) algebras of Monster type.
\begin{flushright} I. Gorshkov \end{flushright}
\end{prob}

The idea behind this question is that if an algebra is $1$-closed, then it should be possible to write explicit product formulas on the space freely generated by axes, and these formulas likely lead to linear relations shrinking this space and hence leading to the final algebra of Monster type. This method has been used successfully for algebras of Jordan type where all algebras are $1$-closed for closed sets of axes.

Related to this is the following question.

\begin{prob}
Is it true that the Frobenius form on an algebra of Monster type exists at least when the algebra is $1$-closed?
\begin{flushright} S. Shpectorov \end{flushright}
\end{prob}

In this case we have a putative form values on a basis consisting of axes (known from the classification of $2$-generated algebras), and we just need to see that this choice of values extends to a full Frobenius form.

\subsection{Double axes and flip subalgebras}

One of the difficulties of studying algebras of Monster type is that we do not know many genuine examples of such algebras outside 
of the case of $2$-generated algebras and also the Griess algebra and its subalgebras, providing examples for the specific 
case $(\al,\bt)=(\frac{1}{4},\frac{1}{32})$. 

Suppose that $M=M_\eta(G,D)$ is a Matsuo algebra constructed from the $3$-transposition group 
$(G,D)$. We refer to the elements of $D$, viewed as elements of $M$, as \emph{Matsuo axes} or \emph{single axes}. 

\begin{defn}
A \emph{double axis} is a sum $a+b\in M$, where $a$ and $b$ are Matsuo axes satisfying $ab=0$. 
\end{defn}

For Matsuo axes the condition $(a,b)=0$ is equivalent to $ab=0$, hence we will call such pairs $a$ and $b$ \emph{orthogonal}.
Double axes are of Monster type $(2\eta,\eta)$ (see \cite{j,gjmss}) provided that $\eta\neq\frac{1}{2}$, but they are not primitive 
in $M$, since $M_1(a+b)=\la a,b\ra$. However, they can be primitive in proper subalgebras of $M$. In particular, they are always 
primitive in the following special subalgebras.

Let $\sg$ be a flip, that is, an automorphism of $(G,D)$ of order $2$. 
Then $\sg$ can be also viewed as an automorphism of $M$. Let $M_\sg=\{u\in M\mid 
u^\sg=u\}$ be the fixed subalgebra. If a double axis $a+b$ is contained in $M_\sg$, with $a,b\notin M_\sg$ then $b=a^\sg$ and $a+b$ is primitive in $M_\sg$.

\begin{defn}
The subalgebra $A(\sg)$ generated by all single and double axes contained in $M_\sg$ is called the \emph{flip subalgebra} corresponding to the flip $\sg$.
\end{defn}

Clearly, flip subalgebras are algebras of Monster type $(2\eta,\eta)$ \cite{j,gjmss} and hence this construction yields a rich class of algebras of Monster type. 
We now turn to the questions of whether algebras of Monster type $(2\eta,\eta)$ can be classified. E.g., can it be that all or nearly all of them come from the flip construction? 
First of all, does the flip construction give us all primitive subalgebras generated by single and double axes?

\begin{prob}
Are there primitive subalgebras of $M$ generated by single and double axes that are not contained in $M_\sg$ for some flip $\sg$? In other words, are there any genuine primitive subalgebras generated by single and double axes that are not obtained by the flip construction?
\begin{flushright} J. McInroy, S. Shpectorov (Question 6.8. \cite{ms})\end{flushright}
\end{prob}

We are aware of some primitive subalgebras of $A(\sg)$ generated by single and double axes that are not equal to the entire $A(\sg)$, so the above question is not about such subalgebras.

Many examples of algebras generated by single and doubles axes have been found and classified by Mamontov and Staroletov \cite{mast}. They covered all such configurations arising from at most five Matsuo axes.

Let $B$ be any primitive subalgebra of $M$ generated by single and double axes.

\begin{defn}
The Matsuo subalgebra generated by all Matsuo axes involved in single and double axes from $B$ is called the \emph{ambient} Matsuo algebra of $B$.
\end{defn}

This is the smallest Matsuo algebra containing $B$.
 
\begin{prob}
Can the ambient Matsuo algebra be recovered from $B$?
\begin{flushright} S. Shpectorov \end{flushright}
\end{prob} 

In other words, can we recover the ambient Matsuo algebra of $B$ without using the embedding of $B$ into $M$?

The concept of the ambient Matsuo algebra leads also to a slight re-evaluation of the flip construction. 
If $B=A(\sg)$ for some flip $\sg$ then the ambient subalgebra $M'$ of $A(\sg)$ is also invariant under $\sg$ and so $\sg$ 
is a flip of $M'$. Furthermore, $A'(\sg)=A(\sg)$. 

\begin{defn}
A flip $\sg$ of a Matsuo algebra $A$ is called \emph{essential} if $A$ is the ambient algebra for $A(\sg)$.
\end{defn}

The concept was introduced in \cite{abs}. Clearly, we should be primarily focusing on essential flips if we 
want to avoid repeatedly constructing the same flip subalgebras.

Hence we pose the following blanket problem.

\begin{prob}
Classify essential flips of $3$-transposition groups.
\begin{flushright} S. Shpectorov \end{flushright}
\end{prob}

We note that the case of irreducible $3$-transposition groups has been mostly completed, although the results still 
need to be properly published. The outstanding irreducible case is the case of the group $G=SO^\eps(n,3)$. The 
reducible case is wide open.

\begin{prob}
Classify (essential) flips and the corresponding flip subalgebras of the Matsuo algebras for groups $SO^\eps(n,3)$.
\begin{flushright} S. Shpectorov \end{flushright}
\end{prob}

The double axis construction can be generalized in an obvious way. Let us say that an idempotent $a$ in an algebra of Jordan type $\eta$ is a $k$-fold axis if $a$ is the sum of $k$ pairwise orthogonal Jordan type axes.

\begin{prob}
What can we say about the fusion law of a $k$-fold axis? Are there any special algebra constructions, similar to the flip construction, that guarantee a nice fusion law and primitivity of $k$-fold axes?
\begin{flushright} C. Franchi, M. Mainardis \end{flushright}
\end{prob}

\subsection{Around the Highwater algebra}

The type $(\al,\bt)=(2,\frac{1}{2})$ is expecially interesting due to presence of the baric and infinite-dimensional Highwater algebra example.
While studying algebras with a binary axes (where all orbits on the axet under the action of the automorphism group have length $2$), McInroy 
and Shpectorov found some non 2-generate examples of baric algebras of Monster type.

\begin{prob}
Can baric algebras of Monster type be classified? Are there infinite-dimensional baric algebras of Monster type $(2,\frac{1}{2})$ for any number of generators?
\begin{flushright} S. Shpectorov \end{flushright}
\end{prob}

Clearly, we are not interested in increasing the number of generators via the direct sum construction.

Let us also mention the following interesting question about the Highwater algebra.

\begin{prob}
Find the identities of the Highwater algebra.
\begin{flushright} I. Kaygorodov\end{flushright}
\end{prob}

\section{Other fusion laws}

So far the focus of research was mostly on the classes of algebras of Jordan and Monster type. However, interesting examples of algebras arise for other fusion laws as well.

\subsection{Natural fusion laws}

Fusion laws are sets with binary operations, not unlike groups and other familiar algebraic structures. What if we re-interpret one of those familiar structures as a fusion law? 

\begin{prob}
The fusion law of an axial algebra is a mapping from the field to the set of all subsets of field. Describe the axial algebras that arise for the extreme fusion law, i.e., the fusion law coming from the multiplication table of the ground field. 
\begin{flushright} I. Gorshkov \end{flushright}
\end{prob}

Clearly, the $0$-eigenspace will be an ideal in the algebra and the $1$-eigenspace is a subalgebra. If the ground field is finite, we will just have a grading by a cyclic group. However, every part of the grading must be a single eigenspace for the specific eigenvalue. Notice however, that the set of elements of the field for which the eigenspaces are non-trivial can be arbitrary.

The following special case looks especially interesting.

\begin{prob}
Describe the axial algebras whose fusion law is the group of roots of unity.   
\begin{flushright} I. Gorshkov \end{flushright}
\end{prob}

In the next question, we do not start with the specific fusion law. Instead, the issue is what the fusion law should look like.

\begin{prob} 
Construct analogues of Matsuo algebras (see Definition \ref{Matsuo}) for arbitrary groups of $p$-transpositions.
\begin{flushright} A. Staroletov \end{flushright}
\end{prob}

For example, consider the dihedral group $G=D_{2p}$, $p$ an odd prime. Similarly to Matsuo algebras, define an algebra with basis $X$, the set of involutions from $G$, and the product, defined for $x,y\in X$, in the following way:
$$
xy=\begin{cases}
x, \text{ if } x=y; \\
\al\cdot{x}+\al\cdot{y}+\sum\limits_{z\in X\setminus\{x,y\}}\bt_z\cdot z,
\end{cases}
$$
for some $\al,\bt_z\in\F$. 

The question above amounts to finding the scalars $\alpha$ and $\beta_z$ that give interesting fusion laws. If $p=5$ then for $\alpha=3\eta$ and $\beta_z=-\eta$ it is easy to see that
$\operatorname{Spec}(ad_s)=\{1,0,4\eta\}$ and eigenspaces of $ad_s$ obey the fusion law in Figure \ref{t:monster}.

\begin{figure}[ht]
\begin{center}
\renewcommand*{\arraystretch}{1.4}
\begin{tabular}{|c||c|c|c|c|}
\hline
$\ast$&$1$&$0$&$4\eta^+$&$4\eta^-$\\
\hline\hline
$1$&$1$&&$4\eta^+$&$4\eta^-$\\
\hline
$0$&&$0$&$4\eta^+$&$4\eta^-$\\
\hline
$4\eta^+$&$4\eta^+$&$4\eta^+$&$1,0$&$4\eta^-$\\
\hline
$4\eta^-$&$4\eta^-$&$4\eta^-$&$4\eta^-$&$1,0,4\eta^+$\\
\hline
\end{tabular}
\end{center}
\caption{Fusion law $\mathcal{M}(4\eta,4\eta)$}\label{t:monster}
\end{figure}

This product also works in the group $5^2:2$, but not in $A_5$ if for products of order $1$, $2$, and $3$ the usual Matsuo multiplication is used.

\subsection{Almost Jordan type fusion law}

Whybrow \cite{w1} did an algorithmic study of all fusion laws of small size. In particular, she observed that for the almost Jordan type fusion law, shown in Table \ref{t:jalbt}, 
\begin{figure}[ht] 
\begin{center}
 \renewcommand{\arraystretch}{1.4}
\begin{tabular}{|c|c|c|c|}
\hline
$\ast$ & 1 & $\al$ & $\bt$ \\ \hline\hline
1 & 1 & $\al$ & $\bt$ \\ \hline
$\al$ & $\al$ & 1, $\al$ & $\bt$ \\ \hline
$\bt$ & $\bt$ & $\bt$ & 1, $\al$ \\ \hline
\end{tabular}
\end{center}
\caption{Fusion law $\cJ(\al,\bt)$ type.}\label{t:jalbt}
\end{figure}
all $2$-generated algebras can be classified. This was developed further by Afanasev \cite{af}, who provided a nearly complete classification of $2$-generated primitive algebras of almost Jordan type $(\al,\bt)$. More recently, Welch \cite{we} completed Afanasiev's classification and he also developed some non 2-generated examples for this fusion law, including generalised Matsuo algebras for $3$-transposition groups.

\begin{prob}
Construct algebras of almost Jordan type for other interesting finite groups. 
\begin{flushright} I. Gorshkov, S. Shpectorov\end{flushright}
\end{prob}

Also, related to Welch's construction of generalised Matsuo algebras, we want to put forward the following question.

\begin{prob}
Is there a construction transforming algebras of Jordan type into algebras of almost Jordan type with $\al\neq 0$?
\begin{flushright} I. Gorshkov, S. Shpectorov\end{flushright}
\end{prob}

Note that among the $2$-generated algebras of almost Jordan type there are algebras not having any analogues among algebras of Jordan type. Hence the construction, as in the question above, cannot be expected to work both ways.

\subsection{Cover of the Highwater algebra}

The cover $\hat\cH$ of the Highwater algebra $\cH$ discovered by Franchi, Mainardis and McInroy \cite{fmm} satisfies the following ``extended Monster type'' fusion law.
\begin{figure}[ht] 
\begin{center}
    \renewcommand{\arraystretch}{1.4}
\begin{tabular}{|c||c|c|c|c|c|}
  \hline
   * & $1$ & $\frac{5}{2}$ & $0$ & $2$ & $\frac{1}{2}$ \\
   \hline \hline
  $1$ & $1$ & $\frac{5}{2}$ & & $2$ & $\frac{1}{2}$ \\
  \hline
  $\frac{5}{2}$ & $\frac{5}{2}$ & $\frac{5}{2}$ & $\frac{5}{2}$ & & $\frac{1}{2}$ \\
  \hline
  $0$ & & $\frac{5}{2}$ & $\frac{5}{2}, 0$ & $\frac{5}{2}, 2$ & $\frac{1}{2}$ \\
   \hline
   $2$ & $2$ & & $\frac{5}{2}, 2$ & $\frac{5}{2}, 0$ & $\frac{1}{2}$ \\
   \hline
   $\frac{1}{2}$ & $\frac{1}{2}$ & $\frac{1}{2}$ & $\frac{1}{2}$ & $\frac{1}{2}$ & $\frac{5}{2}, 0, 2$\\
   \hline
\end{tabular}
\caption{The fusion law $\mathcal{F}$ for $\hat{\mathcal{H}}$}
\end{center}
\end{figure}
Note that in characteristic $5$, the eigenvalue $\frac{5}{2}$ becomes zero, and the fusion law folds into the Monster type $(2,\frac{1}{2})$ fusion law.

\begin{prob}
Is it possible to classify primitive algebras with the above extended Monster type fusion law.
\begin{flushright} C. Franchi, M. Mainardis \end{flushright}
\end{prob}

\subsection{Axial algebras from graphs and geometries}

Cuypers \cite{cu} developed a very interesting construction building axial algebras with fusion laws similar to the Jordan fusion law from graphs and geometries. In particular, he shows that every finite group can be realised as an automorphism group of an axial algebra. In the following question he suggests a concrete way of building axial algebra from graphs or geometries.

\begin{prob}
Consider a geometry (or a graph) and define an algebra having as basis the set of idempotents corresponding to points (vertices) and the product of distinct points $p$ and $q$ being $0$ if they are non-collinear (non-adjacent) and, otherwise, 
$$pq=a(p+q)+b\sum_{r\in\ell} r,$$ 
where $\ell$ is the unique line (or edge) on $p$ and $q$. Here $a$ and $b$ are fixed numbers, parameters of the construction.`
Describe the fusion law of such an algebra.
\begin{flushright} H. Cuypers \end{flushright}
\end{prob}

We note that the geometry is assumed to be a partial linear space, that is, a pair of points is contained in at most one line. In case of Fischer spaces this construction generalises Matsuo algebras. We can also mention that some examples of Welch of algebras of almost Jordan type are as in this question, for specific values of $a$ and $b$.

\medskip%
While the geometry/graph construction is generally quite wild, one could expect that especially nice axial algebras could arise from graphs with additional properties.

\begin{prob}
Are there any specific classes of geometries or graphs leading to good axial algebras. 
\begin{flushright} S. Shpectorov\end{flushright}
\end{prob}

Note that the theory of solid subalgebras and components in algebra has a very strong geometric flavour, where we can view axes as points and $2$-generated subalgebras as lines.

\begin{prob}
Develop a theory of point-line geometries of algebras of Jordan type. 
\begin{flushright} S. Shpectorov\end{flushright}
\end{prob}

\subsection{Axial algebras from analysis and PDEs}

\medskip
Tkachev noticed that all idempotents in Hsiang algebras, arising in the theory of non-linear partial differential equations, satisfy the fusion law shown in Table \ref{t:hsiang}.
\begin{figure}[ht] 
\begin{center}
 \renewcommand{\arraystretch}{1.4}
\begin{tabular}{|c||c|c|c|c|}
\hline
$\ast$ & 1 & $-1$ & $-\frac{1}{2}$ & $\frac{1}{2}$ \\ \hline\hline
1 & 1 & $-1$ & $-\frac{1}{2}$ & $\frac{1}{2}$ \\ \hline
$-1$ & $-1$ & 1 & $\frac{1}{2}$ & $-\frac{1}{2}, \frac{1}{2}$ \\ \hline
$-\frac{1}{2}$ & $-\frac{1}{2}$ & $\frac{1}{2}$ & $1, -\frac{1}{2}$ & $-1, -\frac{1}{2}$ \\ \hline
$\frac{1}{2}$ & $\frac{1}{2}$ & $-\frac{1}{2}, \frac{1}{2}$ & $-1, -\frac{1}{2}$ & $1, -1, -\frac{1}{2}$ \\ \hline
\end{tabular}
\end{center}
\caption{Fusion law $\cH$ for Hsiang algebras}\label{t:hsiang}
\end{figure}

\begin{prob}
Can the $2$-generated algebras with the fusion law $\cH$ be classified?
\begin{flushright} S. Shpectorov\end{flushright}
\end{prob}

The fusion law in Table \ref{t:hsiang} is not $C_2$-graded. At the same time, as Tkachev pointed out, all known Hsiang algebras admit nontrivial automorphism groups. In this respect, the following question might be interesting.

\begin{prob}
Are there any primitive idempotents in Hsiang algebras that satisfy a (stricter) $C_2$-graded sub fusion law of the fusion law $\cH$?
\begin{flushright} S. Shpectorov\end{flushright}
\end{prob}

A possible $C_2$-graded sub fusion law to consider is obtained by removing $-\frac{1}{2}$ from the $2\times 2$ minor in the bottom right corner. This corresponds to $\{1,-1\}$ being the even part of the table and $\{\frac{1}{2},-\frac{1}{2}\}$ the odd part. There could also be other interesting sub laws worth investigating.

\subsection{Chayet-Garibaldi algebras}

It is of course a very interesting question whether all algebraic groups acts on good axial algebras. Desmet \cite{dd} found in Chayet-Garibaldi algebras of type $B_n$ primitive idempotents satisfying the fusion law in Figure \ref{t:bn}.
\begin{figure}[ht] 
\begin{center}
 \renewcommand{\arraystretch}{1.4}
\begin{tabular}{|c||c|c|c|c|}
\hline
$\ast$ & 1 & 0 & $\frac{1}{2}$ & $\frac{1}{4}$ \\ \hline\hline
1 &1 &  & $\frac{1}{2}$ & $\frac{1}{4}$ \\ \hline
0 &  &0  & $\frac{1}{2}, \frac{1}{4}$ & $\frac{1}{2}, \frac{1}{4}$ \\ \hline 
$\frac{1}{2}$ & $\frac{1}{2}$ & $\frac{1}{2}, \frac{1}{4}$  & 1, 0 & 0 \\ \hline
$\frac{1}{4}$ & $\frac{1}{4}$ & $\frac{1}{2}, \frac{1}{4}$  & 0 & 1, 0 \\ \hline
\end{tabular}
\end{center}
\caption{Fusion law $\cB$ for type $B_n$.}\label{t:bn}
\end{figure}

This is a $C_2$-graded fusion law with the partition $\cF_+=\{1,0\}$ and 
$\cF_-=\{\frac{1}{2},\frac{1}{4}\}$. 

\begin{prob}
Is it possible to classify $2$-generated primitive axial algebras with the fusion law $\cB$?
\begin{flushright} S. Shpectorov\end{flushright}
\end{prob}

If this first step is doable then the much more ambitious question would be the following.

\begin{prob}
Are the Chayet-Garibaldi algebras of type $B_n$ characterised by the fusion law $\cB$?
\begin{flushright} S. Shpectorov\end{flushright}
\end{prob}

More generally, all Chayet-Garibaldi algebras, except the one arising for the type $G_2$, have primitive idempotents obeying the fusion law in Figure \ref{t:cnf4} \cite{dd}.
\begin{figure}[ht] 
\begin{center}
 \renewcommand{\arraystretch}{1.4}
\begin{tabular}{|c||c|c|c|c|c|}
\hline
$\ast$ & 1 & 0 & $\frac{1}{2}$ & $\frac{1}{4}$ & $\frac{1}{32}$ \\ \hline\hline
1 &1 &  & $\frac{1}{2}$ & $\frac{1}{4}$ & $\frac{1}{32}$ \\ \hline
0 &  & 0 & $\frac{1}{2}, \frac{1}{4}, \frac{1}{32}$ & $\frac{1}{2}, \frac{1}{4}$ & $\frac{1}{2}, \frac{1}{32}$ \\ \hline

$\frac{1}{2}$ & $\frac{1}{2}$ & $\frac{1}{2}, \frac{1}{4}, \frac{1}{32}$  & 1, $\frac{1}{4}, \frac{1}{32}, 0$ & $\frac{1}{2}, \frac{1}{32}, 0$ & $\frac{1}{2}, \frac{1}{4}, 0$ \\ \hline

$\frac{1}{4}$ & $\frac{1}{4}$ & $\frac{1}{2}, \frac{1}{4}$  & $\frac{1}{2}, \frac{1}{32}, 0$ & 1, 0 & $\frac{1}{2}, \frac{1}{32}$ \\ \hline
$\frac{1}{32}$ & $\frac{1}{32}$ & $\frac{1}{2}, \frac{1}{32}$  & $\frac{1}{2}, \frac{1}{4}, 0$ & $\frac{1}{2}, \frac{1}{32}$ & 1, $\frac{1}{4}, 0$ \\ \hline
\end{tabular}
\end{center}
\caption{Fusion law $\cC\cG$ for type $C_n$, $F_4$.}\label{t:cnf4}
\end{figure}
Existence of primitive idempotents satisfying $\cC\cG$ in Chayet-Garibaldi algebras of type $C_n$ and $F_4$ was shown in \cite{dd}. The fusion law $\cB$ for type $B_n$ is a sub law of $\cC\cG$, so this case is also covered. 

The type $G_2$ is a true exception. In this case, \cite{dd} cites the fusion law in Figure \ref{t:g2} arising for certain primitive idempotents.
\begin{figure}[ht] 
\begin{center}
    \renewcommand{\arraystretch}{1.4}
\begin{tabular}{|c||c|c|c|c|}
\hline
$\ast$ & 1 & 0 & $\frac{1}{2}$ & $\frac{1}{3}$ \\ \hline\hline
1 & $1$ &  & $\frac{1}{2}$ & $\frac{1}{3}$ \\ \hline
0 &  & $0$ & $\frac{1}{2}, \frac{1}{3}$ & $\frac{1}{2}, \frac{1}{3}$ \\ \hline

$\frac{1}{2}$ &  $\frac{1}{2}$ & $\frac{1}{2}, \frac{1}{3}$ & $1, 0, \frac{1}{3}$ & $\frac{1}{2}, 0$ \\ \hline
$\frac{1}{3}$ & $\frac{1}{3}$ & $\frac{1}{2}, \frac{1}{3}$ & $\frac{1}{2}, 0$ & $1, 0, \frac{1}{3}$ \\ \hline
\end{tabular}
\end{center}
\caption{Fusion law $\cG$ for type $G_2$.}\label{t:g2}
\end{figure}
Interestingly, this fusion law is not $C_2$-graded.

We will formulate the following blanket questions concerning these fusion laws.

\begin{prob}
What can we say about algebras with fusion laws $\cC\cG$ and $\cG$? E.g., can $2$-generated primitive algebras with these fusion laws be classified?
\begin{flushright} S. Shpectorov\end{flushright}
\end{prob}

Note that the fusion law $\cC\cG$ is amazingly close to the Monster fusion law 
$\cM(\frac{1}{4},\frac{1}{32})$. Namely, the Monster fusion law is a sub fusion law of $\cC\cG$ obtained by removing the eigenvalue $\frac{1}{2}$. Note that 
$\frac{1}{2}$ is responsible for having an automorphism group of positive dimension (equivalently, having nontrivial algebra Lie of derivations). Hence, from this point of view, the axial algebras of Monster type, including the Griess algebra, can be viewed as finite analogues of the Chayet-Garibaldi algebras. 

In this respect the following question seems very interesting.

\begin{prob}
Construct $2$-graded axial algebras whose Miyamoto group is isomorphic to a pariah sporadic simple group. 
\begin{flushright} E. Zelmanov \end{flushright}
\end{prob}

Of particular interest is the question of whether there exist such algebras of Monster type or, perhaps, with the fusion law $\cC\cG$.

\section{Computational aspects}

Almost simple groups of $6$-transpositions have recently been classified by Parker and Turner \cite{pt}.

\begin{prob}
Determine for which almost simple $6$-transposition groups $G$ there exists an algebra of Monster type with some parameters $\al$, $\bt$. 
\begin{flushright} S. Shpectorov \end{flushright}
\end{prob}

The list of Parker and Turner is quite long and a comprehensive solution for the above problem may not be feasible. However, any new examples of this kind would certainly be very interesting. 

In fact, this question is also interesting for the stricter Jordan type fusion law. Let us call a (finite) group a \emph{Jordan type group} if it arises as the Miyamoto group for some algebra of Jordan type half. In view of the available results on components in algebras of Jordan type half, for an almost simple Miyamoto group $G$, it is almost certain that the underlying algebra $A$ is a Jordan algebra and hence the finite axet leading to $G$ is a small part of the 
infinite set of all primitive axes of $A$. Also, in general a finite group of Jordan type does not need to be a $6$-transposition groups.

Still the following problem is very interesting. 

\begin{prob}
Compute algebras of Jordan type $\frac{1}{2}$ for known $6$-transposition groups. 
\begin{flushright} S. Shpectorov \end{flushright}
\end{prob}

The key example, found computationally by McInroy and Shpectorov, is the $15$-dimensional Jordan algebra over $\Q(\sqrt{5})$ on which the sporadic Hall-Janko simple group $J_2$
arises as a Miyamoto group of an axet of size $315$, in a bijection with the set of involutions in $J_2$. Note that $J_2$ is a $6$-transposition group.

Staroletov noticed that if $G$ arises on an algebra of Jordan type half over $\Q$ then it is a $6$-transposition group missing product order $5$. The list of Parker and Turner contains a number of such examples. 

Other smaller examples found by McInroy and Shpectorov are all groups of real and complex groups of reflections. It is almost certain that the existence of the $J_2$ example has something to do with the fact that $J_2$ is a group of quaternionic reflections. For real reflections groups $G$ with a simply laced diagram and the underlying space $V$, De Medts and Rehren \cite{dr} showed that the Jordan algebra of symmetric bilinear forms on $V$ arises as a factor algebra of the Matsuo algebra of $G$. (Note that $G$ in this case is a $3$-transposition group.) 

The following question is theoretical in nature, but it is strongly aligned with the above computational programme.

\begin{prob}
Is it true that the construction of De Medts and Rehren extends to complex and quaternioning groups of reflections?
\begin{flushright} S. Shpectorov \end{flushright}
\end{prob}

If the answer to this question is positive then we can ask whether all examples of Jordan type groups arise in these way or whether there are Jordan type groups beyond reflection groups.

Finally, there is a group of questions related to Ivanov's programme of studying Majorana algebras, which are a special case of algebras of Monster type $(\frac{1}{4},\frac{1}{32})$. Many of these algebras arise as subalgebras of the Griess algebras. 

\begin{prob}
Construct the Griess algebra explicitly in a computer and use the computer representation to determine the subalgebras corresponding to smaller sporadic groups and other interesting subgroups of Monster.
\begin{flushright} S. Shpectorov \end{flushright}
\end{prob}

Note that the standard realisation of the Griess algebra, via the tensor product, is likely beyond what the modern computers can do. However, it may be possibly to realise the Griess algebra product directly, say, as described by Conway \cite{c}, 
even though this still remains a very defficult problem.

As an intermediate project, one can try to construct explicitly an interesting subalgebra of the Griess algebra.

\begin{prob}
Construct explicitly the algebra for $(A_{12}\times A_5):2$, or at least the algebra for $S_{12}$, which is a subalgebra of the above.
\begin{flushright} S. Shpectorov \end{flushright}
\end{prob}

The subgroup $(A_{12}\times A_5):2$ is the top of the chain of subgroups featuring subgroups $A_n$, $n\leq 12$, and their centralisers in the Monster. The subalgebras for the subgroups $A_n$ from this series were studied extensively by Franchi, Ivanov and Mainardis \cite{fim2,fim,fm}, however even the dimensions of the larger of these subalgebras are not firmly known. An explicit construction, based on the wealth of information collected in \cite{fim2,fim}, will help to confirm all dimensions and it can also help with further related problems.

Looking at the known and conjectured dimensions of the algebras from the $A_n$ series, one may notice that growth in dimension observed for small $n$, slows down as $n$ approaches $12$. Conjecturally, this may be due to the fact that positive definiteness of the Frobenius form is one of the axioms of Majorana algebras and this condition is used in \cite{fim}. Can it be that there is a larger algebra for $S_{12}$, whose Frobenius form is positive semi-definite?

\begin{prob}
Does the algebra for $S_{12}$ have a larger cover?
\begin{flushright} S. Shpectorov \end{flushright}
\end{prob}

The radical of such an algebra would be an ideal and the factor over it would be the algebra $S_{12}$ that we know as a subalgebra of the Griess algebra.

Currently the series of these $A_n$ algebras stops at $n=12$. If the above cover were to exist, this would open the possibility of algebra for groups $A_n$ with $n>12$ with indefinite Frobenius form.

Note that a known factor algebra puts strong restrictions on the structure of a possible cover algebra. In particular, it would be interesting to incorporate this special case in the expansion algorithm \cite{ms2}, which was used to compute algebras for many smaller groups. 

\begin{prob}
Realise a version of the expansion algorithm to compute covers of known algebras.
\begin{flushright} S. Shpectorov \end{flushright}
\end{prob}

The information from the known factor can feed into the expansion algorithm to yield additional relations and hence compute the larger algebras in this way.

More generally, while the expansion algorithm was very successful in computing many algebras for known groups, it would be good to develop it further, so it can handle larger groups and algebras. One of the ideas in this respect is to use the module structure of the algebra, which would allow to group relations in sizable chunks, instead of dealing with them individually.

\begin{prob}
Realise the version of the expansion algorithm computing relations module-wise.
\begin{flushright} S. Shpectorov \end{flushright}
\end{prob}

The expansion algorithm works by expanding the known part of the axial algebra using multi-linear algebra operations such the tensor product and symmetric square of modules. These can be precomputed for the irreducible modules that are potentially involved in the algebra and then, these can be used as building blocks.

The interesting case where these ideas can be employed is the Monster type algebras for symmetric and alternating groups, where there is a very good description of all irreducible modules.

Another related idea for improving the expansion algorithm is to work over a function field where indeterminates represent the relevant parameters of the algebra. This can be used to extract additional relations earlier, which should speed up the algorithm. In particular, this idea can be used for splitting the $1$- and $0$-eigenvectors and utilising the primitivity assumption.

\begin{prob}
Realise the version of the expansion algorithm over function fields.
\begin{flushright} S. Shpectorov \end{flushright}
\end{prob}

Note that the current MAGMA version of the expansion algorithm due to McInroy does not make any assumptions about the ground field and so, in principle, it will work over the function field. So the problem is more about utilising the additional parameters within the algorithm in order to obtain additional relations.


Whybrow \cite{w1} created a GAP code for computing $2$-generated axial algebras with a given small fusion law, where the elements of the fusion law can be viewed as parameters. 

\begin{prob}
Restore Wybrow's code (and possibly port it into MAGMA) and continue the exploration of small fusion laws.
\begin{flushright} S. Shpectorov \end{flushright}
\end{prob}
\section{Acknowledgements}
	The work of the first author was supported the Russian Science Foundation, project 24-11-00119, https://rscf.ru/project/24-11-00119/ .

\Addresses
 
\end{document}